\documentclass[12pt]{article}
\usepackage[pdftex, plainpages=false]{hyperref}
\usepackage{amsfonts,amsmath}
\usepackage{graphicx}
\def \beq {\begin{eqnarray}}
\def \eeq {\end{eqnarray}}
\def \beqn {\begin{eqnarray*}}
\def \eeqn {\end{eqnarray*}}

\newcommand{\halmos}{\rule{1ex}{1.4ex}}
\newcounter{for}[section]

\newtheorem{itlemma}{Lemma}[section]
\newtheorem{itproposition}[itlemma]{Proposition}
\newtheorem{theorem}[itlemma]{Theorem}
\newtheorem{itcorollary}[itlemma]{Corollary}
\newtheorem{itremark}[itlemma]{Remark}
\newtheorem{itremarks}[itlemma]{Remarks}
\newtheorem{itdefinition}[itlemma]{Definition}
\newtheorem{itexample}[itlemma]{Example}

\newenvironment{fact}{\begin{itfact}\rm}{\end{itfact}}
\newenvironment{claim}{\begin{itclaim}\rm}{\end{itclaim}}
\newenvironment{lemma}{\begin{itlemma}}{\end{itlemma}}
\newenvironment{remark}{\begin{itremark}\it}{\end{itremark}}
\newenvironment{remarks}{\begin{itremarks}\rm}{\end{itremarks}}
\newenvironment{corollary}{\begin{itcorollary}}{\end{itcorollary}}
\newenvironment{proposition}{\begin{itproposition}}{\end{itproposition}}
\newenvironment{definition}{\begin{itdefinition}}{\end{itdefinition}}
\newenvironment{example}{\begin{itexample}\rm}{\end{itexample}}
\newenvironment{proof}{\noindent {\em Proof}.\ \
}{\hspace*{\fill}$\halmos$\medskip}
\newcommand{\be}[1]{\addtocounter{for}{1} \begin{equation}\label{#1}}
\newcommand{\ee}{\end{equation}}
\newcommand{\bl}[1]{\begin{lemma}\label{#1}}
\newcommand{\br}[1]{\begin{remark}\label{#1}}
\newcommand{\brs}[1]{\begin{remarks}\label{#1}}
\newcommand{\bt}[1]{\begin{theorem}\label{#1}}
\newcommand{\bd}[1]{\begin{definition}\label{#1}}
\newcommand{\bp}[1]{\begin{proposition}\label{#1}}
\newcommand{\bc}[1]{\begin{corollary}\label{#1}}
\newcommand{\bfact}[1]{\begin{fact}\label{#1}}
\newcommand{\bex}[1]{\begin{example}\label{#1}}
\newcommand{\ec}{\end{corollary}}
\newcommand{\efact}{\end{fact}}
\newcommand{\eex}{\end{example}}
\newcommand{\el}{\end{lemma}}
\newcommand{\er}{\end{remark}}
\newcommand{\ers}{\end{remarks}}
\newcommand{\et}{\end{theorem}}
\newcommand{\ed}{\end{definition}}
\newcommand{\ep}{\end{proposition}}
\newcommand{\epr}{\end{proof}}
\newcommand{\bpr}{\begin{proof}}
\newcommand{\bcl}[1]{\begin{claim}\label{#1}}
\newcommand{\ecl}{\end{claim}}
\newcommand{\ecs}{\end{corollary}}
\newcommand{\eers}{\end{exercise}}
\newcommand{\eexs}{\end{example}}
\newcommand{\eems}{\end{example}}
\newcommand{\els}{\end{lemma}}
\newcommand{\eles}{\end{lemmaex}}
\newcommand{\ets}{\end{theorem}}
\newcommand{\eds}{\end{definition}}
\newcommand{\eps}{\end{proposition}}
\newcommand{\bi}{\begin{itemize}}
\newcommand{\ei}{\end{itemize}}
\newcommand{\ben}{\begin{enumerate}}
\newcommand{\een}{\end{enumerate}}

\def\vbar{\mathchoice{\vrule height6.3ptdepth-.5ptwidth.8pt\kern-.8pt}
   {\vrule height6.3ptdepth-.5ptwidth.8pt\kern-.8pt}
   {\vrule height4.1ptdepth-.35ptwidth.6pt\kern-.6pt}
   {\vrule height3.1ptdepth-.25ptwidth.5pt\kern-.5pt}}
\def\fudge{\mathchoice{}{}{\mkern.5mu}{\mkern.8mu}}
\def\bbc#1#2{{\rm \mkern#2mu\vbar\mkern-#2mu#1}}
\def\bbb#1{{\rm I\mkern-3.5mu #1}}
\def\bba#1#2{{\rm #1\mkern-#2mu\fudge #1}}
\def\bb#1{{\count4=`#1 \advance\count4by-64 \ifcase\count4\or\bba A{11.5}\or
   \bbb B\or\bbc C{5}\or\bbb D\or\bbb E\or\bbb F \or\bbc G{5}\or\bbb H\or
   \bbb I\or\bbc J{3}\or\bbb K\or\bbb L \or\bbb M\or\bbb N\or\bbc O{5} \or
   \bbb P\or\bbc Q{5}\or\bbb R\or\bbc S{4.2}\or\bba T{10.5}\or\bbc U{5}\or
   \bba V{12}\or\bba W{16.5}\or\bba X{11}\or\bba Y{11.7}\or\bba Z{7.5}\fi}}

\def \qed {{\hspace*{\fill}$\halmos$\medskip}}

\def \T {{\cal{T}}}
\def \TT{{\mathsf{T}}}
\def \L {{\cal{L}}}

\def \A {{\cal{A}}}

\def \G {{\cal{G}}}
\def \E {{\cal{E}}}
\def \C {{\cal{C}}}
\def \D {{\cal{D}}}
\def \Db{\underline{{\mathcal D}}_{b-\epsilon}}
\def \Dh{{\bar{\mathcal D}}_{b+\epsilon}}

\def \S {{\cal{S}}}

\newcommand{\ba}[1]{\addtocounter{for}{1} \begin{eqnarray}\label{#1}}
\newcommand{\ea}{\end{eqnarray}}

\def\sqr#1#2{{\vcenter{\vbox{\hrule height .#2pt
                             \hbox{\vrule width .#2pt height#1pt \kern#1pt
                                   \vrule width .#2pt}
                             \hrule height .#2pt}}}}

\def\pmb#1{\setbox0=\hbox{#1}%
   \kern-.025em\copy0\kern-\wd0
   \kern.05em\copy0\kern-\wd0
 \kern-.025em\raise.0433em\box0 }
\def\sqr#1#2{{\vcenter{\vbox{\hrule height.#2pt
     \hbox{\vrule width.#2pt height#1pt \kern#1pt
   \vrule width.#2pt}\hrule height.#2pt}}}}

\def\B{{\mathcal B}}
\def \RR{{\mathcal R}}

\def \I{{\mathcal I}}
\def\N{{\mathbb N}}
\def\Z{{\mathbb Z}}
\def\R{{\mathbb R}}
\def\Q{{\mathbb Q}}

\def\P{{\mathbb P}}

\def\EE{{\mathbb E}}

\def\diam{\text{diam}}
\def\dist{\text{dist}}
\def\capa{\text{cap}}
\def\Weight{\text{Weight}}

\def\e{\epsilon}

\def\v{\varphi}

\def\bs{\backslash}
\def\reff#1{(\ref{#1})}

\def \ind {\hbox{1\hskip -3pt I}}
\newcommand {\cro}[1] {\left[ {#1} \right]}
\newcommand {\acc}[1] {\left\{ {#1} \right\}}
\newcommand {\pare}[1] {\left( {#1} \right)}

\newcommand {\bra}[1] {\left< {#1} \right>}

\begin{document}

\title{Large Deviation Principle for Self-Intersection Local Times 
for Random Walk in $\Z^d$ with $d\ge 5$.}
\author{Amine Asselah \\ Universit\'e Paris-Est\\
amine.asselah@univ-paris12.fr}
\date{}
\maketitle
\begin{abstract}
We obtain a large deviation principle for the self-intersection
local times for a symmetric random walk in dimension $d\ge 5$. As an application,
we obtain moderate deviations for random walk in random sceneries in Region II
of \cite{AC05}.
\end{abstract}



{\em Running head}: LDP for self-intersections in $d\ge 5$.

\section{Introduction.} \label{sec-intro}
We consider an aperiodic symmetric random walk on the lattice $\Z^d$, with
$d\ge 5$. More precisely, if $S_n$ is the position of the 
walk at time $n\in \N$, then
$S_{n+1}$ chooses uniformly at random a site of $\acc{z\in \Z^d: |z-S_n|\le 1}$,
where for $z=(z_1,\dots,z_d)\in \Z^d$, the $l^1$-norm is $|z|:=|z_1|+\dots+|z_d|$.
When $S_0=x$, we denote the law of this walk by $P_x$, and its expectation by $E_x$.

We are concerned with estimating the number of trajectories of length $n$
with {\it many} self-intersections, in the large $n$-regime. The
self-intersection local times process reads as follows
\be{def-SILT}
\text{for } n\in \N,\qquad B_n=\sum_{0\le i<j< n} \ind\{ S_i=S_j\}.
\ee
The study of self-intersection local times has a long history in probability
theory, as well as in statistical physics. Indeed, a caricature of a polymer would
be a random walk self-interacting
through short-range forces; a simple model arises as we penalize the 
simple random walk law with
$\exp(\beta B_n)$, where $\beta<0$ corresponds to a weakly self-avoiding walk,
and $\beta>0$ corresponds to a self-attracting walk. The question is whether there is
a transition from collapsed paths to diffusive paths,
as we change the parameter $\beta$. 
We refer to Bolthausen's Saint-Flour
notes~\cite{SF-bolt} for references and a discussion of these models.

It is useful to represent $B_n$ in terms of local times
$\acc{l_{n}(x), x\in \Z^d}$, 
that is the collection of number of visits of $x$ up to time $n$,
as $x$ spans $\Z^d$. We set, for $k<n$,
\be{intro.1}
l_{[k,n[}(x)=\ind\{S_k=x\}+\dots+\ind\{S_{n-1}=x\},
l_n=l_{[0,n[},\quad\text{and}\quad 
||l_n||_2^2=\sum_{z\in \Z^d} l^2_n(z).
\ee
It is immediate that $||l_n||_2^2=2 B_n+n$. Henceforth, we always
consider $||l_n||_2^2$ rather than $B_n$. It turns out useful
to think of the self-intersection local times
as the square of the $l^2$-norm of an additive and positive process
(see Section~\ref{sec-ergodic}).
Besides, we will deal with other $q$-norm of $l_n$
(see Proposition~\ref{prop-alpha}),
for which there is no counterpart in terms of multiple self-intersections.

In dimensions $d\ge 3$, a random walk spends, 
on the average, a time of the order
of one on most visited sites, whose number, up to time $n$, is of order $n$.
More precisely, a result of \cite{BS} states
\be{intro.3}
\frac{1}{n}||l_n||_2^2\stackrel{L^2}{\longrightarrow}
\gamma_d=2G_d(0)-1,
\quad\text{with}\quad \forall z\in \Z^d,\ G_d(z)=\sum_{n\ge 0} P_0(S_n=z).
\ee
The next question concerns estimating the probabilities of large deviations 
from the mean: that is $P_0(||l_n||_2^2-E_0[||l_n||_2^2]\ge n\xi)$
with $\xi>0$. In dimension $d\ge 5$, the speed
of the large deviations is $\sqrt n$, and we know from
\cite{AC05} that a finite (random) set of sites, say $\D_n$, visited of the order of
$\sqrt n$ makes a dominant contribution to produce the excess self-intersection.

However, in dimension 3, 
the correct speed for our large deviations is $n^{1/3}$ (see \cite{A06}), and
the excess self-intersection is made up by sites visited less than
some power of $\log(n)$. It is expected that the walk spends most of its time-period
$[0,n]$ on a ball of radius of order $n^{1/3}$. 
Thus, in this box, sites are visited a time of order unity.

The situation is still different in dimension 2. First,
$E_0[B_n]$ is of order $n\log(n)$, and a result of Le Gall~\cite{legall} states that
$\frac{1}{n}(B_n-E_0[B_n])$ converges in law to a non-gaussian random variable.  
The large (and moderate) deviations asymptotics obtained recently
by Bass, Chen \& Rosen in~\cite{bcr05b}, reads as follows.
There is some positive constant $C_{BCR}$, such that for any sequence
$\acc{b_n,n\in \N}$ going to infinity 
with $\lim_{n\to\infty} \frac{b_n}{n}=0$, we have
\be{asympt-bcr}
\lim_{n\to\infty} \frac{1}{b_n} \log\pare{P\pare{ B_n-E_0[B_n]\ge b_n n}}=-C_{BCR}.
\ee
For a LDP in the case of $d=1$, we refer to Chen and Li~\cite{chen-li}
(see also Mansmann~\cite{Ma91} for the case of a Brownian motion instead of a random walk).
In both $d=2$ and $d=1$, the result is obtained by showing that the local times
of the random walk is close to its smoothened conterpart.

Finally, we recall a related result of Chen and M\"orters \cite{chen-morters}
concerning mutual intersection local times of two independent random walks
in infinite time horizon when $d\ge 5$. Let $l_\infty(z)=\lim_{n\to\infty}
l_n(z)$, and denote by $\tilde l_\infty$ an independent copy of $l_\infty$.
All symbols related to the second walk differ with a tilda.
We denote the average over both walks by $\EE$, and the product law is denoted
$\P$.
The intersection local times of two random walks, in an infinite time horizon, is
\[
\bra{l_\infty,\tilde l_\infty}=\sum_{z\in \Z^d} l_\infty(z)\tilde l_\infty(z),\quad
\text{and}\quad \EE\cro{\bra{l_\infty,\tilde l_\infty}}=\sum_{z\in \Z^d}
G_d(z)^2<\infty,
\]
where Green's function, $G_d$, is square summable in dimension 5 or more.
Chen and M\"orters in~\cite{chen-morters}
have obtained sharp asymptotics for $\{\langle l_\infty,\tilde l_\infty
\rangle\ge t\}$ for $t$ large, in dimension 5 or
more, by an elegant asymptotic estimation of the moments, improving on the
pioneering work of Khanin, Mazel, Shlosman and Sinai in~\cite{KMSS}. Their
method provides a variational formula for the rate functional, and
their proof produces (and relies on) a finite volume version.
Namely, for any finite subset $\Lambda\subset
\Z^d$,
\be{CM-finite}
\lim_{t\to\infty} \frac{1}{\sqrt t} \log \P
\pare{ \bra{\ind_{\Lambda} l_\infty,\tilde l_\infty} \ge t}
=-2\I_{CM}(\Lambda),\quad\text{and}\quad \lim_{\Lambda \nearrow \Z^d}
\I_{CM}(\Lambda)=\I_{CM},
\ee
with
\[
\I_{CM}=\inf\acc{ ||h||_2:\ h\ge 0,\ ||h||_2<\infty, \text{and}\ ||U_h||\ge 1},
\]
where
\be{rate-CM}
U_h(f)(x)=\sqrt{e^{h(x)}-1}\sum_{y\in \Z^d} \pare{G_d(x-y)-\delta_x(y)}
(f(y) 
\sqrt{e^{h(y)}-1}),
\ee
and $\delta_x$ is Kronecker's delta function at $x$.

In this paper, we consider self-intersection local-times, and
we establish a Large Deviations Principle in $d\ge 5$. 
\bt{intro-th.1} We assume $d\ge 5$. There is a constant $\I(2)>0$, such that for $\xi>0$
\be{intro.4}
\lim_{n\to\infty} 
\frac{1}{\sqrt{ n}}\log 
P_0\pare{||l_n||_2^2-E\cro{||l_n||_2^2}\ge n \xi}=-\I(2) \sqrt{\xi}.
\ee
Moreover,
\be{identification}
\I(2)=\I_{CM}.
\ee
\et
\br{rem-2statements}
The reason for dividing Theorem~\ref{intro-th.1} into two statements
\reff{intro.4} and \reff{identification} is that our proof has two steps:
(i) The proof of the existence of the limit in \reff{intro.4}, which
relies eventually on a subadditive argument, in spite of an odd scaling;
(ii) An identification with the constant of Chen and M\"orters.

Also, we establish later the existence of a limit for other $q$-norms 
of the local-times (see Proposition~\ref{prop-alpha}), 
for which we have no variational formulas.
\er

The identification \reff{identification} relies on the fact that
both the excess self-intersection local times and 
large intersection local times
are essentially realized on a finite region. This is explained
heuristically in Remark 1 of \cite{chen-morters}, and we provide
the following mathematical statement of this latter phenomenon.
\bp{prop-CM} Assume dimension is 5 or more.
\be{ineq-CM}
\limsup_{\epsilon\to 0} \limsup_{t\to\infty} \frac{1}{\sqrt t}\log
\P\pare{\sum_{z\in \Z^d}
\ind_{\{\min(l_\infty(z),\tilde l_\infty(z))<\epsilon \sqrt t\}}
l_\infty(z)\tilde l_\infty(z)>t}=-\infty.
\ee
\ep

%Application to other powers.
% Application to RWRS.
Finally, we present applications of our results to 
Random Walk in Random Sceneries (RWRS). 
We first describe RWRS. We consider a field $\{\eta(x),x\in \Z^d\}$
independent of the random walk $\{S_k,k\in \N\}$, and made up of
symmetric unimodal i.i.d.\  with law denoted by $\Q$ and tail decay characterized by
an exponent $\alpha>1$ and a constant $c_{\alpha}$ with
\be{eq-tail.1}
\lim_{t\to\infty} \frac{\log \Q\pare{\eta(0)>t}}{t^{\alpha}}=-c_{\alpha}.
\ee
The RWRS is the process
\[
\bra{\eta,l_n}:= \sum_{z\in \Z^d} \eta(z) l_n(z)= \eta(S_0)+\dots+
\eta(S_{n-1}).
\]
We refer to \cite{AC05} for references for RWRS, and for a diagram of the speed of
moderate deviations $\acc{\bra{\eta,l_n}> \xi n^\beta}$ with $\xi>0$,
in terms of $\alpha>1$ and $\beta>\frac{1}{2}$. In this paper,
we concentrate on what has been called in \cite{AC05} Region II:
\be{def-RegII}
1< \alpha<\frac{d}{2},\quad\text{and}\quad
1-\frac{1}{\alpha+2}<\beta<1+\frac{1}{\alpha}.
\ee
In region II, the random walk is expected to visit often a few sites, and it is therefore
natural that our LDP allows for better asymptotics in this
regime. We set
\be{old-defIIbis}
\zeta=\beta\frac{\alpha}{\alpha+1}(<1),\qquad \frac{1}{\alpha^*}=1-\frac{1}{\alpha},
\quad\text{and for $\chi>0$}\quad \bar\D_n(\xi):=\acc{z:\ l_n(z)\ge \xi}.
\ee
In bounding from above the probability of $\acc{ \bra{\eta,l_n}\ge \xi\ n^{\beta}}$, we take
exponential moments of $\bra{\eta,l_n}$, and first integrate 
with respect to the $\eta$-variables.
Thus, the behavior of the log-Laplace transform of $\eta$, say $\Gamma(x)=\log
E\cro{\exp(x\eta(0))}$, either at zero or at infinity, plays a key r\^ole. 
This, in turn,
explains why we need a LDP for other powers of the
local times. For $q\ge 1$, the $q$-norm of function $\v:\Z^d\to\R$ is
\[
||\v||_{q}^q:=\sum_{z\in \Z^d} |\v(z)|^{q}.
\]
Before dealing with $\acc{\bra{\eta,l_n}> \xi n^\beta}$, we give estimates for 
the ${\alpha^*}$-norm of the local-times, for $\alpha^*>\frac{d}{d-2}$.
\bp{prop-alpha} 
Choose $\zeta$ as in \reff{old-defIIbis} with $\alpha,\beta$ in Region II.
Choose $\chi$ such that  $\zeta> \chi\ge \frac{\zeta}{d/2}$, and any $\xi>0$.
There is a positive constant $\I(\alpha^*)$ such that
\be{alpha-star}
\lim_{n\to\infty} \frac{1}{n^\zeta}
\log\pare{ P\pare{ ||\ind_{\bar\D_n(n^\chi)} l_n||_{\alpha^*} 
\ge  \xi n^\zeta}}=- \xi\ \I(\alpha^*).
\ee
\ep
Our moderate deviations estimates for RWRS is as follows.
\bt{prop-rwrs}
Assume  $\alpha,\beta$ are in Region II given in \reff{def-RegII}. With $\zeta$
given in \reff{old-defIIbis}, and any $\xi>0$
\be{rwrs-md}
\lim_{n\to\infty} \frac{1}{n^\zeta} \log\pare{
P\pare{\bra{\eta,l_n}\ge  \xi  n^\beta}}=
-c_{\alpha} (\alpha+1) \pare{ \frac{\I(\alpha^*) }{\alpha}}^{\frac{\alpha}{\alpha+1}}
\quad \xi^{\frac{\alpha}{\alpha+1}}.
\ee
\et

We now wish to outline schematically the main ideas and limitations
in our approach. This serves also to describe the organisation of the
paper. First, we use a 
shorthand notation for the centered self-intersection local times process,
\be{def-silt}
\overline{||l_n||_2^2}=||l_n||_2^2-E_0\cro{||l_n||_2^2}.
\ee
Theorem~\ref{intro-th.1} relies on the following intermediary
result interesting on its own.
\bp{prop-sub}
Assume $d\ge 5$. There is $\beta>0$, such that for any $\epsilon>0$,
there is $\alpha_\epsilon>0$, and $\Lambda_\epsilon$
a finite subset of $\Z^d$, such that for any $\alpha>\alpha_\epsilon$,
for any $\Lambda\supset \Lambda_\epsilon$ finite, and $n$ large enough
\be{lawler.10}
\begin{split}
\frac{1}{2} P_0\big(&||\ind_\Lambda l_{\lfloor\alpha\sqrt{n}\rfloor}
||_2^2\ge n\xi (1+\epsilon),\ S_{\lfloor\alpha\sqrt{n}\rfloor}=0\big)\\
&\le P_0(\overline{||l_n||_2^2}\ge n\xi) \le e^{\beta \epsilon \sqrt n}
P_0\pare{||\ind_\Lambda l_{\lfloor\alpha\sqrt{n}\rfloor}
||_2^2\ge n\xi (1-\epsilon),\ S_{\lfloor\alpha\sqrt{n}\rfloor}=0}.
\end{split}
\ee
We use $\lfloor x\rfloor$ for the integer part of $x$.
\ep
The upper bound for $P_0(\overline{||l_n||_2^2}\ge n\xi)$
in \reff{lawler.10} is the main technical result of the paper.

From our previous work in \cite{AC05}, 
we know that the main contribution
to the excess self-intersection comes from level set 
$\D_n=\{x:l_{n}(x)\sim {\sqrt n}\}$. This is the place
where $d\ge 5$ is crucial. Indeed, this latter fact is false in dimension
3 as shown in \cite{A06}, and unknown in $d=4$. 
In Section~\ref{sec-old}, we recall and refine the results of \cite{AC05}.
We establish that $\D_n$ is a {\it finite} set.
More precisely, for any $\e>0$ and $L$ large
enough, there is a constant $C_\e$ such that for $n$ large enough
\be{strat.1}
P\pare{\overline{||l_n||_2^2}\ge n\xi}\le C_\e\ P\pare{
||\ind_{\D_n}l_n||_2^2\ge n\ \xi(1-\epsilon),\ |\D_n|<L}.
\ee
Then, our main objective is to show that the time spent
on $\D_n$ is of order $\sqrt n$. However, this is only
possible if some control on the diameter of $\D_n$ is first
established. This is the main difficulty.
Note that $\D_n$ is visited by the random 
walk within the time-period $[0,n[$,
and from \reff{strat.1}, a crude uniform estimate yields
\be{strat.2}
P\pare{\overline{||l_n||_2^2}\ge n\xi}
\le C_\e (2n)^{dL}
\sup_{\Lambda\in ]-n,n[^d,|\Lambda|\le L}\!\!\!
 P\pare{ ||\ind_\Lambda l_n||_2^2\ge n\ \xi(1-\epsilon)}.
\ee
Now, we can replace the time period $[0,n[$, in the right hand side
of \reff{strat.2}, by an infinite interval $[0,\infty)$
since the local time increases with time.
Consider $\Lambda_n\subset ]-n,n[^d$ which realizes
the supremum in \reff{strat.2}. 
Next, we construct two maps: a {\it local} map 
$\T$ in Section~\ref{sec-moveC}, and a {\it global}
map $f$ in Section~\ref{sec-moveloop}.
A finite number of
iterates of $\T$ (at most $L$), say $\T^L$, transforms $\Lambda_n$ 
into a subset of finite diameter. On the other hand,
$f$ maps $\{\D_n=\Lambda_n\}$ into $\{\D_n=\T(\Lambda_n)\}$, allowing
us to compare the probabilities of these two events.
Thus, the heart of our argument has two ingredients.
\begin{itemize}
\item A {\it marriage theorem} which is recalled in 
Section~\ref{sec-marriage}. 
It is then used to perform {\it global surgery} on the circuits.
\item Classical potential estimates of Sections~\ref{sec-cageloop} and
\ref{sec-movetrip}. This is the place where the random 
walk's features enter the play. Our estimates relies on
basic estimates (Green's function asymptotics,
Harnack's inequalities and heat kernel asymptotics), which
are known to hold for general symmetric random walks (see
\cite{lawler-limic}). Though we have considered the simplest
aperiodic symmetric random walk, all our results hold when the basic
potential estimates hold.
\end{itemize}
We then iterate $f$ a finite number of time to reach
$\{\D_n=\T^L(\Lambda_n)\}$. 
To control the cost of this transformation, it is crucial that
only a finite number of iterations of $f$ is needed. 
The construction of $\T$ and $f$ requires as well many preliminary
steps.
\begin{enumerate}
\item Section~\ref{sec-cluster} deals with {\it clusters}.
In Section~\ref{sec-cluster}, we introduce a partition of $\Lambda_n$ 
into a collection of nearby points, called {\it clusters}. 
In Section~\ref{sec-moveC}, we define a map
$\T$ acting on {\it clusters}, by translating one {\it cluster} at a time.
\item Section~\ref{sec-circuit} deals with {\it circuits}.
In Section~\ref{sec-defcir},
we decompose a trajectory in $\{\D_n=\Lambda_n\}$
into all possible {\it circuits}. We introduce the notions of
{\it trip} and {\it loop}.
\end{enumerate}
We show in Proposition~\ref{prop-time.1},
that for trajectories in $\{\D_n=\T^L(\Lambda_n)\}$,
no time is wasted on lengthy excursions, and
the total time needed to visit $\D_n$ 
is less than $\alpha{\sqrt n}$, for some large $\alpha$. This
steps also relies  on assuming $d\ge 5$. Indeed, we have been using that
conditionned on returning to the origin, the expected return time
is finite in dimension 5 or more.
This concludes the outline of the proof of 
the upper bound in Proposition~\ref{prop-sub}. The lower bound is
easy, and is done in Section~\ref{sec-LB}.

Assuming Proposition~\ref{prop-sub}, 
we are in a situation where a certain $l^2$-norm of an additive process
is larger than $\sqrt{n\xi}$ over a time-period of $\alpha {\sqrt n}$.
Section~\ref{sec-subadditive} presents
a subadditive argument yielding the existence of a limit \reff{intro.4}.
We identify the limit in Section~\ref{sec-ident}. We prove
Proposition~\ref{prop-CM} in Section~\ref{sec-KMSS}.
Finally, the proof of Theorem~\ref{prop-rwrs} is given 
in Section~\ref{sec-rwrs}.

We conclude by
mentionning two outstanding problems out of our reach.
\begin{itemize}
\item Establish a Large Deviations Principle in $d=3$, showing that the walk spends most
of its time during time-period $[0,n[$, in a ball of radius about $n^{1/3}$. 
\item In dimension 4, find which level set of the local times gives a dominant contribution
to making the self-intersection large.
\end{itemize}

\section{Preliminaries on Level Sets.}\label{sec-old}
In this section, we recall and refine the analysis of \cite{AC05}. 
The approach of \cite{AC04,AC05} focuses on the contribution of each
level set of the local times to the 
event $\{||l_n||_2^2-E[||l_n||_2^2]>n \xi\}$. 
This section is essentially a corollary of \cite{AC05}. 

We first recall Proposition 1.6 of \cite{AC05}. For $\epsilon_0>0$, set
\[
\RR_n=\acc{x\in \Z^d: n^{1/2-\epsilon_0}\le l_n(x)\le n^{1/2+\epsilon_0}}.
\]
Then, for any $\epsilon>0$
\be{18-AC05}
\lim_{n \rightarrow \infty} \frac{1}{\sqrt{n}} \log
P \pare{ ||\ind_{\RR_n^c} l_n||_2^2-E_0\cro{||l_n||_2^2}
\geq n\epsilon\xi} = - \infty \, .
\end{equation}
Thus, we have for any $0<\epsilon< 1$, and $\xi>0$
\be{level-main1}
P\pare{\overline{||l_n||_2^2}\geq n \xi}\le 
P \pare{ ||\ind_{\RR_n^c} l_n||_2^2-E_0\cro{||l_n||_2^2}
\geq n\epsilon\xi}+
P\pare{ ||\ind_{\RR_n} l_n||_2^2 \ge n\xi(1-\epsilon)}.
\ee
We only need to focus on the second term of the right hand side 
of \reff{level-main1},
and for simplicity here, we use $\xi>0$ instead of $\xi(1-\epsilon)$. 
First, we show in Lemma~\ref{level-lem.1} that when asking
$\{||l_n||_2^2 \geq E[||l_n||_2^2] +n\xi\}$ with $\xi>0$, we can assume
$\{||l_n||_2^2 \le An\}$ for some large $A$.
Then, in Lemma~\ref{level-lem.2}, 
we show that the only sites which matter are those 
whose local times is within $[\frac{\sqrt{n}}{A}, A\sqrt{n}]$ for some large constant $A$. 

\bl{level-lem.1}
For $A$ positive, there are constants $C,\kappa>0$ such that
\be{level.1}
P\pare{\overline{||l_n||_2^2}\ge nA}\le C \exp\pare{-\kappa \sqrt{An}}.
\ee
\el
\bpr
We rely on Proposition 1.6 of \cite{AC05}, and the proof of Lemma 3.1 of \cite{AC05}
(with $p=2$ and $\gamma=1$), 
for the same subdivision $\acc{b_i,i=1,\dots,M}$ of $[1/2-\epsilon,1/2]$, 
and the same $\acc{y_i}$ such that $\sum y_i\le 1$,
but the level sets are here of the form
\be{level.2}
\D_i=\acc{x\in \Z^d:\ A^{\frac{1}{2}} n^{b_i}\le l_n(x)< A^{\frac{1}{2}} 
n^{b_{i+1}}}.
\ee
Using Lemma 2.2 of \cite{AC05}, we obtain the second line of \reff{level.3},
\ba{level.3}
P\pare{ \sum_{\cup \D_i} l_n^2(x) \ge nA}&\le&\sum_{i=1}^{M-1} P\pare{
|\D_i|(A^{\frac{1}{2}}n^{b_{i+1}})^2\ge ny_i A}=
\sum_{i=1}^{M-1} P\pare{|\D_i|\ge y_i n^{1-2b_{i+1}}}\cr
&\le & \sum_{i=1}^{M-1} (n^d)^{n^{1-2b_{i+1}}y_i} 
\exp\pare{ -\kappa_d A^{\frac{1}{2}} n^{b_{i}+(1-\frac{2}{d})
(1-2b_{i+1})}y_i^{1-\frac{2}{d}}}\cr
&\le & \sup_{i\le M}\acc{\C_i(n) 
\exp\pare{ -\kappa_d A^{\frac{1}{2}} n^{b_{i}+(1-\frac{2}{d})
(1-2b_{i+1})}y_i^{1-\frac{2}{d}}}},
\ea
where $\C_i(n):=M (n^d)^{n^{1-2b_{i+1}}y_i}$.
The constant $\kappa_d$ is linked with estimating the probability of spending a given time
in a given domain $\Lambda$ of prescribed volume; this latter inequality 
is derived in Lemma 1.2 of \cite{AC04}.
We first need $\C_i(n)$ to be negligible, which imposes
\be{level.4}
n^{1-2b_{i+1}}y_i \log(n^d)\ll A^{\frac{1}{2}} 
n^{b_{i}+(1-\frac{2}{d})(1-2b_{i+1})}y_i^{1- \frac{2}{d}}
\ee
Inequality \reff{level.4} is easily seen to hold when $b_i$ is larger than $1/2-\epsilon$,
for $\epsilon$ small. Now, we need that for some $\kappa>0$
\be{level.5}
\kappa_dA^{\frac{1}{2}}
n^{b_{i}+(1-\frac{2}{d})(1-2b_{i+1})}(y_i)^{1- \frac{2}{d}}
\ge 2\kappa A^{\frac{1}{2}}\sqrt{n}.
\ee
This holds with the choice of $y_i$ as in Lemma 3.1 of \cite{AC05}.
We use one $\kappa$ of \reff{level.5} to match 
$\C_i(n)$ in \reff{level.3}, and we are left with a constant $C$ such that
\be{level.7}
P\pare{ \sum_{\cup \D_i} l_n^2(x) \ge nA}\le C \exp\pare{-\kappa \sqrt{An}}
\ee
\epr

% NEED TO WRITE HERE THE LEMMA WHICH SPEAKS ABOUT FINITE SIZE SQUARE ROOT OF N
For any positive reals $A$ and $\zeta$, an $k\in \N\cup\{\infty\}$, we define
\be{level.15}
\D_k(A,\xi):=\acc{x\in \Z^d: \frac{\xi}{A}\le l_k(x)< A\xi}.
\ee
\bl{level-lem.2} Fix $\xi>0$. For any $M>0$, there is $A>0$ so that 
\be{level.16}
\limsup_{n\to\infty} \frac{1}{\sqrt n} \log\pare{P\pare{ \sum_{\RR_n\bs\D_n(A,\sqrt n)}
l_n^2(x)>n\xi}}\le -M.
\ee
Also, 
\be{level.17}
P\pare{ | \D_n(A,\sqrt n)|\ge A^3}\le C\exp\pare{-\kappa \sqrt{ An}}.
\ee
\el
\bpr
We consider an increasing sequence $\acc{a_i,i=1,\dots,N}$ to be chosen later, and form
\be{level.8}
\B_i=\acc{x:\ \frac{\sqrt n}{a_i}\le l_n(x)< \frac{\sqrt n}{a_{i-1}}},
\ee
where $a_0$ will be chosen as a large constant, and $a_N\sim n^{\epsilon}$.
In view of Lemma~\ref{level-lem.1}, it is enough to show
that the probability of the event $\acc{\sum_{\B_i}
l_n^2(x)\ge n\xi}$ is negligible. First, from Lemma~\ref{level-lem.1}, we can 
restrict attention to $\acc{An\ge \sum_{\B_i} l_n^2(x)\ge n\xi_i}$ for some 
large constant $A$ and with $\xi=\sum \xi_i$ a decomposition to be chosen later. 
When considering the sum over $x\in  \B_i$, we obtain
\be{level.9}
\sum_{x\in \B_i} l_n^2(x)\le nA\Longrightarrow |\B_i|\pare{\frac{\sqrt n}{a_i}}^2\le An
\Longrightarrow |\B_i|\le a_i^2 A.
\ee
Similarly, we obtain the 
lower bound $|\B_i|\ge \xi_i a_{i-1}^2$. If we call
\be{level.10}
H_i=\acc{  a_{i-1}^2\xi_i<|\B_i|\le a_i^2 A},
\ee
then by Lemma~\ref{level-lem.1}, if we set $l_n(\B_i)=\sum_{x\in \B_i} l_n(x)$
\ba{level.11}
P\pare{ \sum_{x\in \B_i} l_n^2(x)> n\xi_i}
&\le &P(\sum_{x\in \B_i} l_n^2(x)> n A))+P\pare{H_i\cap \acc{l_n(\B_i)\ge a_{i-1}\xi_i
\sqrt{n}}}\cr
&\le & C e^{-\kappa \sqrt{An}}
+ (n^d)^{a_i^2 A} \exp\pare{- \kappa_d \frac{ a_{i-1}\xi_i\sqrt{n}}{(a_i^2A)^{2/d}}}.
\ea
Since we assume $a_i\le n^{\epsilon}$, 
the term $ (n^d)^{a_i^2 A}$ is innocuous. It remains
to find, for any large constant $M$, 
two sequences $\acc{a_i,\xi_i,i=1,\dots,N}$ such that
\be{level.12}
\kappa_d \frac{ a_{i-1}\xi_i}{(a_i^2A)^{2/d}}=M,\quad\text{and}\quad
\sum \xi_i=\xi.
\ee
Fix an arbitrary $\delta>0$ and set
\be{level.13}
a_i:=(1+\delta)^i a_0,
\quad
\xi_i:= \frac{z(\delta)}{(1+\delta)^{\gamma i}} \xi,
\quad\text{and}\quad \gamma=1-\frac{4}{d},
\ee
where $z(\delta)$ is a normalizing constant ensuring that $\sum \xi_i=\xi$.
Using the values \reff{level.13} in \reff{level.12}, we obtain
\be{level.14}
\frac{\kappa_dz(\delta) \xi}{(1+\delta) A^{2/d}} a_0^{1-4/d} =M.
\ee
Now, for any constant $M$, we can choose an $a_0$ large enough so that none of the
level $\B_i$ contributes. Note also that $N=\min\acc{n: a_n\ge n^\epsilon}$.

Finally, \reff{level.17} follows from Lemma~\ref{level-lem.1}, once we note that 
\[
P\pare{|\D_n(A,\sqrt n)|\ge A^3}\le 
P\pare{||\ind_{\D_n(A,\sqrt n)} l_n||_2^2\ge An}.
\]
\epr

We will need estimates for other powers of the local times. 
We choose two parameters $(\alpha,\beta)$ satisfying \reff{def-RegII}, and
we further define
\be{old-defII}
\zeta=\beta\frac{\alpha}{\alpha+1},\qquad b=\frac{\beta}{\alpha+1},\qquad
\frac{1}{\alpha^*}=1-\frac{1}{\alpha},
\quad\text{and}\quad \bar\D_n(n^b):=\acc{z:\ l_n(z)\ge n^b}.
\ee

When dealing with the $\alpha^*$-norm of $l_n$, we only focus on sites with large
local times. Among those sites, we show that finitely many contribute to
making the $\alpha^*$-norm of $l_n$ {\it large}. To appreciate the first estimate,
similar in spirit and proof to Lemma~\ref{level-lem.2}, 
recall that $\zeta<1$, $\alpha^*>1$, and
$||l_n||_{\alpha^*}\ge ||l_n||_1=n$.
\bl{level-lem.3} 
Choose $\zeta,b$ as in \reff{old-defII} with $\alpha,\beta$ in Region II.
For any $\xi>0$, there are constants $C,\kappa>0$ such that
\be{level.22}
P\pare{ ||\ind_{\bar\D_n(n^b)} 
l_n||_{\alpha^*} \ge \xi n^\zeta }\le C \exp\pare{-\kappa \xi n^{\zeta}}.
\ee
Moreover, for any $M>0$, there is $A>0$ such that
\be{level.20}
\limsup_{n\to\infty} \frac{1}{n^\zeta} \log P\pare{ ||\ind_{\D_n(A,n^\zeta)^c\cap 
\bar\D_n(n^b)} l_n||_{\alpha^*}>\xi n^{\zeta}}\le -M.
\ee
Finally, from \reff{level.22}, we have
\be{level.21}
P\pare{ | \D_n(A,n^\zeta)|\ge A^2}\le 
P\pare{ ||\ind_{\D_n(A,n^\zeta)} l_n||_{\alpha^*}\ge An^\zeta}\le
C\exp\pare{-\kappa An^{\zeta}}.
\ee
\el
The proof is similar to that of Lemmas~\ref{level-lem.1} and ~\ref{level-lem.2}, and we omit 
the details. We point out that Lemma 3.1 of \cite{AC05} has to be used with
$p=\alpha^*$ and $\gamma=\alpha^*\zeta$. Also, Proposition 3.3 of \cite{AC05} holds on
$\bar\D_n(n^b)$ since the condition $b \frac{d}{2}\ge \zeta$ is fulfilled in Region II.

\section{Clusters' Decomposition.}\label{sec-cluster}
From Lemma~\ref{level-lem.2}, for any $\epsilon>0$, and
$A$ large enough, we have $C_\epsilon>0$ such that for 
any $\xi>0$ and $n$ large enough
\be{circuit.1}
P\pare{\overline{||l_n||_2^2}\ge n\xi(1+\epsilon)} \le 
C_\epsilon P\pare{||\ind_{\D_n(A,\sqrt n)} l_n||_2^2\ge n 
\xi,\ |\D_n(A,\sqrt n)|\le A^3}.
\ee
Since $\D_n(A,\sqrt n)\subset ]-n,n[^d$,
we bound the right hand side of \reff{circuit.1} by a uniform bound
\be{circuit.2}
\begin{split}
P\pare{\overline{||l_n||_2^2}\ge n\xi(1+\epsilon)} \le&
C_\e (2n)^{dA^3} \sup_{\Lambda}
P\pare{||\ind_{\Lambda} l_n||_2^2\ge n \xi,\D_n(A,\sqrt n)=\Lambda}\\
\le&
C_\e (2n)^{dA^3} \sup_{\Lambda}
P\pare{||\ind_{\Lambda} l_{\infty}||_2^2\ge n \xi,
\Lambda\subset \D_{\infty}(A,\sqrt n)},
\end{split}
\ee
where in the supremum over $\Lambda$ we assumed
that $\Lambda\subset ]-n,n[^d,\ |\Lambda|\le A^3$. Also,
in $\D_{\infty}(A,\sqrt n)$ (defined in \reff{level.15})
we may adjust with a larger $A$ if necessary.

If we denote by $\Lambda_n$ the finite subset of $\Z^d$ which
realizes the last supremum in \reff{circuit.2},
then our starting point, in this section, is
the collection $\{\Lambda_n,n\in \N\}$ of finite subsets of $\Z^d$.

\subsection{Defining Clusters.}\label{sec-defC}
In this section, we partition an arbitrary finite
subset of $\Z^d$, say $\Lambda$ into subsets of nearby sites, 
with the feature that these subsets are far apart. More precisely,
this partitioning goes as follows.

\bl{lem-cluster} Fix $\Lambda$ finite subset of $\Z^d$, and $L$ an
integer. There is a partition of $\Lambda$ whose elements are called
$L$-clusters with the property that
two distinct $L$-clusters $\C$ and $\tilde \C$ satisfy
\be{clus.3}
\dist(\C,\tilde \C):=\inf\acc{|x-y|,\ x\in\C,y\in \tilde \C}\ge
4\max\pare{ \diam(\C), \diam(\tilde\C),L}.
\ee
Also, there is a positive constant $C(\Lambda)$ which
depends on $|\Lambda|$, such that for any $L$-cluster $\C$
\be{clus.key2}
\diam(\C)\le C(\Lambda)\ L.
\ee
\el
\br{rem-cluster} If we define an $L$-shell $\S_L(\C)$ around $\C$ by
\be{clus.key1}
\S_L(\C)=\acc{z\in \Z^d:
\ \dist(z,\C)\le \max(L,\diam(\C))},\quad\text{ then }
\quad \S_L(\C)\cap \Lambda=\C. 
\ee
We deduce from \reff{clus.key2},
and \reff{clus.key1}, 
that for any $\C$ and any
$x,y\in \C$, there is a finite sequence of points $x_0=x,\dots,x_k=y$ 
(not necessarely in $\Lambda$), such that for $i=1,\dots,k$
\be{clus.key3}
|x_i-x_{i-1}|\le L,\quad\text{and}\quad B(x_i,L)\subset \S(\C)\quad
(\text{where }B(x_i,L)=\acc{z\in \Z^d:\ |x_i-z|\le L}).
\ee
\er
\bpr
We build clusters by a bootstrap algorithm.
At level 0, we define a {\it linking} relation for $x,y\in \Lambda$:
$x\overset{0}{\leftrightarrow} y$ if $|x-y|\le 4L$, and an equivalent relation
$x\overset{0}{\sim} y$ if there is a (finite) path $x=x_1,x_2,\dots,x_k=y\in \Lambda$ such
that for $i=1,\dots,k-1$, $x_i\overset{0}{\leftrightarrow} x_{i+1}$.
The cluster at level 0 are the equivalent classes of $\Lambda$. We denote
by $\C^{(0)}(x)$ the class which contains $x$, and by $|\C^{(0)}|$ the number
of clusters at level 0 which is bounded by $|\Lambda|$. It is important to note that the diameter
of a cluster is bounded independently of $n$. Indeed, it is easy to see,
by induction on $|\Lambda|$, that for any $x\in \Lambda$, we have
$\diam(\C^{(0)}(x))\le 4L(|\C^{(0)}(x)|-1)$, so that 
\be{clus.4}
\diam(\C^{(0)}(x))\le 4L |\Lambda|.
\ee
Then, we set
\be{clus.1}
x\overset{1}{\leftrightarrow} y\quad\text{if}\quad |x-y|\le 4 \max\pare{ \diam(\C^{(0)}(x)),
\diam(\C^{(0)}(y)),L}.
\ee
As before, relation $\overset{1}{\leftrightarrow}$ is associated with an equivalence
relation $\overset{1}{\sim}$ which defines clusters $\C^{(1)}$. Note also that 
$x\overset{0}{\sim} y$ implies that $x\overset{1}{\sim}y$, and that for any $x\in \Lambda$, 
\be{clus.diam}
\diam(\C^{(1)}(x))\le 5 |\C^{(0)}|
\max\acc{ \diam(\C):\ \C\in \C^{(0)}}\le 5 |\Lambda| (4L |\Lambda|),
\ee
since we produce $\C^{(1)}$'s by multiple concatenations of 
pairs of $\C^{(0)}$-clusters at a distance
of at most four times the maximum diameters of the clusters making up level 0, 
those latter clusters being less in number than $|\Lambda|$. 
In the worst scenario, there is
one cluster at level 1 made up of all clusters of $\C^{(0)}$ at a distance of at most 
$4 \max\acc{\diam(\C):\ \C\in \C^{(0)}}$.

If the number of clusters
at level 0 is the same as those of level 1, then the algorithm stops and
we have two distinct clusters $\C,\tilde \C\in \C^{(0)}$
\[
\dist(\C,\tilde \C):=\inf\acc{|x-y|,\ x\in\C,y\in \tilde \C}\ge
4\max\pare{ \diam(\C), \diam(\tilde\C),L}.
\]
Otherwise, the number
of cluster at level 1 has decreased by at least one. Now, assume by way of induction,
that we have reached level $k-1$. We define $\overset{k}{\leftrightarrow}$ as follows
\be{clus.2}
x\overset{k}{\leftrightarrow} y\quad\text{if}\quad |x-y|\le 4 \max\pare{ \diam(\C^{(k-1)}(x)),
\diam(\C^{(k-1)}(y)),L}.
\ee
Now, since $|\Lambda|$ is finite, the 
algorithm stops in a finite number of steps.
The clusters we obtain eventually are called $L$-clusters. 
Note that two distinct $L$-clusters satisfy \reff{clus.3}. Property 
\reff{clus.key2} with $C(\Lambda)=(5|\Lambda|)^{|\Lambda|}$, 
follows by induction with the same argument used to prove \reff{clus.diam}.
\epr

\subsection{Transforming Clusters.}\label{sec-moveC}
For a subset $\Lambda$ and an integer $L$, assume that we have
a partition in terms of $L$-cluster as in Lemma~\ref{lem-cluster}.
We define the following map on the partition of $\Lambda$.
\bl{lem-move} There is a map $\T$ on the $L$-clusters
of $\Lambda$ such that $\T(\C)=\C$, but for one cluster, say $\C_1$
where $\T(\C_1)$ is a translate of $\C_1$
such that, when the following minimun is taken over all $L$-clusters
\be{red-cluster}
0=\min\acc{\dist\pare{\C,\T(\C_1)}-\pare{\diam(\C)+\diam(\T(\C_1)}}.
\ee
Also, for any $L$-cluster $\C\not= \C_1$, we have
\be{notbad-cluster}
\dist(\C,\T(\C_1))\le 2 \dist(\C,\C_1).
\ee
\el
We denote by $\T(\Lambda)=\cup\  \T(\C)$. Also, we can define
$\T$ as a map on $\Z^d$: for a site $z\in \C_1$
$\T(z)$ denotes the translation of $z$, otherwise 
$\T(z)=z$. Finally, we can define the inverse of $\T$, which we
denote $\T^{-1}$.
\br{move-rem.1}
Note that $\T(\Lambda)$ has at least one $L$-cluster less than $\Lambda$ 
since \reff{clus.3} does not hold for $(\C_0,\T(\C_1))$. 
Thus, if we apply to $L$-cluster partition of Lemma~\ref{lem-cluster}
to $\T(\Lambda)$, $\C_0$ and $\T(\C_1)$ would
merge into one $L$-cluster, possibly triggering other merging.
\er
\bpr
We start with two clusters which minimize
the distance among clusters. Let $\C_0$ and $\C_1$ be such that
\be{clus.5}
\dist(\C_0,\C_1)= \min\acc{\dist(\C,\C'):\ \C,\C' 
\text{ distinct clusters}}.
\ee
Now, let $(x_0,x_1)\in \C_0\times\C_1$ such that 
$|x_0-x_1|=\dist(\C_0,\C_1)$, and note
that by \reff{clus.3}, $|x_0-x_1|\ge 2\pare{ \diam(\C_0)+\diam(\C_1)}$.
Assume that $\diam(\C_0)\ge \diam(\C_1)$. 
We translate sites of $\C_1$ by a vector
whose coordinates are the integer parts of the following vector
\be{clus.6}
u=(x_0-x_1)\pare{1-\frac{\diam(\C_0)+\diam(\C_1)}{|x_0-x_1|}},
\ee
in such a way that the translated cluster, say $\T(\C_1)$, is at a distance
$\diam(\C_0)+\diam(\C_1)$ of $\C_0$.
We now see that $\T(\C_1)$ is far enough from other clusters.
Let, as before, $z\in \C$, and note that
\ba{clus.8}
|z-\tilde y|&\ge& |z-x_0|-|x_0-\tilde y|\ge |z-x_0|-\pare{|x_0-\tilde x_1|+|\tilde x_1-\tilde y|}\cr
&\ge & 4 \max(\diam(\C),\diam(\C_0))-\pare{ \diam(\C_0)+\diam(\C_1)+\diam(\C_1)}\cr
&\ge & \diam(\C)+\diam(\C_1)
\ea
Thus, for any cluster $\C$, we have 
\be{clus-diam}
\dist(\C,\T(\C_1))\ge \diam(\C)+\diam(\T(\C_1)).
\ee
Finally, we prove \reff{notbad-cluster}.
Let $z$ belong to 
say $\C\not= \C_1$, and let $\tilde y\in \T(\C_1)$ be the image
of $y\in \C_1$ after translation by $u$. 
Then, using that $\dist(\C_0,\C_1)$ minimizes
the distance among distinct clusters
\ba{clus.7}
|z-\tilde y|&\le& |z-y|+|y-\tilde y|\le |z-y|+\dist(\C_0,\C_1)\cr
&\le& |z-y|+\dist(\C,\C_1)\le 2|z-y|.
\ea
\epr
\section{ On Circuits.}\label{sec-circuit}
\label{circuit}
\subsection{Definitions and Notations.}\label{sec-defcir}
Let $\Lambda_n\subset\Z^d$ maximizes the supremum in the last term
of \reff{circuit.2}.
Assume we have partitioned $\Lambda_n$ into $L$-clusters,
 as done in Section~\ref{sec-cluster}.

We decompose the paths realizing $\{||1_{\Lambda_n
}l_\infty||_2^2\ge n\xi\}$ 
with $\{\Lambda_n\subset\D_\infty(A,\sqrt n)\}$ 
into the successive visits to $\Lambda_n'=\Lambda_n\cup \T(\Lambda_n)$. 
For ease of notations, we drop the subscript $n$ in $\Lambda$ though
it is important to keep in mind that $\Lambda$ varies as we increase
$n$.

We consider the collection of integer-valued vectors over
$\Lambda'$ which we think of as candidates 
for the local times over $\Lambda'$. Thus
\be{circuit.3}
V(\Lambda',n):=\acc{ {\bf k}\in \N^{\Lambda'}:\ 
\inf_{x\in \Lambda}k(x)\ge \frac{\sqrt n}{A},\ \sup_{x\in \Lambda'}k(x)\le
A{\sqrt n},\quad \sum_{x\in \Lambda} k^2(x)\ge n\xi}.
\ee
Also, for ${\bf k}\in V(\Lambda',n)$, we set 
\be{main-estimate}
|{\bf k}|=\sum_{x\in \Lambda'} k(x),\quad\text{and note that}\quad
|{\bf k}|\le |\Lambda'|A{\sqrt n}\le 2A^4 {\sqrt n}.
\ee

We need now more notations. For 
$U\subset \Z^d$, we call $T(U)$ the first hitting time of $U$, and we denote
by $T:=T(\Lambda')=\inf\acc{n\ge 0: S_n\in \Lambda'}$. We also use the notation
$\tilde T(U)=\inf\acc{n\ge 1: S_n\in U}$.
For a trajectory in the event $\acc{l_\infty(x)=k(x),\forall x\in \Lambda'}$, we call 
$\acc{T^{(i)},i\in \N}$ the successive times of visits of $\Lambda'$: 
$T^{(1)}=\inf\acc{n\ge 0: S_n\in \Lambda'}$, and by induction for $i\le |{\bf k}|$ when
$\acc{T^{(i-1)}<\infty}$
\be{circuit.4}
T^{(i)}=\inf\acc{n>T^{(i-1)}: S_n\in \Lambda'}.
\ee
The first observation is that the number of {\it long trips} cannot be too large.
\bl{circuit-lem.1} For any $\epsilon>0$, and $M>0$, there is $L>0$ such that 
for each ${\bf k}\in V(\Lambda',n)$,
\be{circuit.5}
P\pare{ l_\infty|_{\Lambda'}
={\bf k},\ \big|\acc{i\le |{\bf k}|:\ |S_{T^{(i)}}-S_{T^{(i-1)}}\big|>
{\sqrt L}}|\ge \epsilon {\sqrt n}}\le e^{-M{\sqrt n}}.
\ee
\el
We know from \cite{AC05} that the probability that 
$\{\overline{||l_n||_2^2}\ge n\xi\}$ is bounded from below
by $\exp(-\bar c \sqrt n)$ for some positive constant $\bar c$. 
We assume $M>2\bar c$ (and $L>L(M)$ given in Lemma~\ref{circuit-lem.1}),
and the left hand side of \reff{circuit.5} is negligible.
The proof of this Lemma is postponed to the Appendix.

We consider now the collections of possible sequence of visited 
sites of $\Lambda'$, and in view of Lemma~\ref{circuit-lem.1}, 
we consider at most $\epsilon {\sqrt n}$ 
consecutive sites at a distance larger than ${\sqrt L}$. 
First, for ${\bf k}\in V(\Lambda',n)$, and each
${\bf z}\in \E({\bf k})$, and $x\in \Z^d$, 
we denote by $l_{{\bf z}}(x)$ the
{\it local times} of ${\bf z}$ at $x$, that is the number of 
occurrences of $x$ in the string ${\bf z}$. Then,
\be{circuit.6}
\E({\bf k})=\acc{{\bf z}\in (\Lambda')^{|{\bf k}|}:\ l_{{\bf z}}(x)=k(x),
\forall x\in \Lambda',\quad \sum_{i< |{\bf k}|} 
\ind_{\acc{|z(i+1)-z(i)|>{\sqrt L}}}< \epsilon {\sqrt n}}.
\ee
\bd{def-circuit}
For ${\bf k}\in V(\Lambda',n)$, a {\it circuit} is an element of
$\E({\bf k})$. The random walk follows 
circuit ${\bf z}\in \E({\bf k})$, if it belongs to the event
\be{circuit.13}
\acc{S_{T^{(i)}}=z(i),i=1,\dots,|{\bf k}|}\cap \acc{ 
T^{(|{\bf k}|+1)}=\infty}.
\ee
\ed
When we lift the second constrain in \reff{circuit.13}, 
we obtain when $L$ is large enough (with the convention $z(0)=0$)
\be{circuit.14}
P\pare{||\ind_{\Lambda}l_\infty||_2^2\ge n\xi,\ \Lambda
\subset \D_\infty(A,\sqrt n)}\le
2\sum_{{\bf k}\in V(\Lambda',n)}
\sum_{{\bf z}\in\E({\bf k})} \prod_{i=1}^{|{\bf k}|}
P_{z(i-1)}\pare{S_T=z(i)}.
\ee
We come now to the definitions of {\it trips} and {\it loops}.
\bd{def-trip} Let ${\bf k}\in V(\Lambda',n)$ and ${\bf z}\in \E({\bf k})$. 
A {\it trip} is a pair
$(z(i),z(i+1))$, where $z(i)$ and $z(i+1)$ 
do not belong to the same cluster. A {\it loop} is a maximal substring
of ${\bf z}$ belonging to the same cluster.
\ed
\br{rem-trip}
We think of a circuit as a succession of loops connected by trips. 
Recall that \reff{clus.3} tells us that two points of a trip
are at a distance larger than $L$. Thus, trips
are necessarely long journeys, whereas loops may contain many 
short journeys, typically of the order of $\sqrt n$. 
For ${\bf z}\in \E({\bf k})$, the number of trips is less 
than $\epsilon {\sqrt n}$, so is the number of loops, since a loop is
followed by a trip.
\er
We recall the notations of
Section~\ref{sec-moveC}: $\Lambda=\{\C_0,\C_1,\dots,\C_k\}$ with
$\dist(\C_0,\C_1)$ minimizing distance
among the clusters. The map $\T$ translates
only cluster $\C_1$. 

%For simplicity, we relabel $\C=\C_1$ and $\tilde \C=\T(\C_1)$. 
%
We now fix ${\bf k}\in V(\Lambda',n)$ and ${\bf z}\in \E({\bf k})$. 
We number the different points of entering and exiting from $\C_1$. 
\be{circuit.7}
\tau_1=\inf\acc{n>0:\ z(n)\in \C_1},\quad\text{and}\quad
\sigma_1=\inf\acc{n>\tau_1:\ z(n)\not\in \C_1},
\ee
and by induction, if we assume $\acc{\tau_2,\sigma_2,\dots,\tau_i,\sigma_i}$ defined
with $\sigma_i<\infty$, then
\be{circuit.8}
\tau_{i+1}=\inf\acc{n>\sigma_i:\ z(n)\in \C_1},\quad\text{and}\quad
\sigma_{i+1}=\inf\acc{n>\tau_{i+1}:\ z(n)\not\in \C_1}.
\ee

\bd{def-loop} For a configuration ${\bf z}\in \E({\bf k})$, its
$i$-th $\C_1$-loop is
\be{circuit.9}
\L(i)=\acc{z(\tau_i),z(\tau_i+1),\dots,z(\sigma_i-1)}.
\ee
We associate with $\L(i)$ the entering and exiting site from $\C_1$, 
$p(i)=\acc{z(\tau_i),z(\sigma_i-1)}$, 
which we think of as the {\it type} of the $\C_1$-loop.
\ed

The construction is identical for $\T(\C_1)$
(usually with a tilda put on all symbols).

\subsection{Encaging Loops.}\label{sec-cageloop}
We wish eventually to transform a piece of random walk
associated with a $\C_1$-loop, into a piece of random walk associated with
a $\T(\C_1)$-loop. We explain one obvious problem we face when
acting with $\T$ on circuits. Consider a $\C_1$-loop in a circuit ${\bf z}$.
Assume for simplicity, that it corresponds to the $i$-th $\C_1$-loop.
In general, 
\be{T-problem}
\prod_{k=\tau_i}^{\sigma_i-2} P_{z(k)}\pare{S_T=z(k+1)}\quad\not=\quad
\prod_{k=\tau_i}^{\sigma_i-2} P_{\T(z(k))}\pare{S_T=\T(z(k+1))}.
\ee
However, if while travelling from $z(k)$ to $z(k+1)$,
the walk were forced to stay inside an
$L$-shell of $\C_1$ during $[\tau_i,\sigma_i[$, then
under $\T$, we would have a walk travelling
from $\T(z(k))$ to $\T(z(k+1))$, inside an $L$-shell of $\T(\C_1)$.

To give a precise meaning to our use of the expression {\it encage},
we recall that for any cluster $\C$, the $L$-shell around $\C$ is denoted
\[
\S(\C)=\acc{z:\ \dist(z,\C)=\max(L,\diam(\C))}.
\]
Now, for $x,y\in \C$, the random walk is encaged inside $S$ while
flying from $x$ to $y$ if it does not exit $\S$ before touching $y$.
The main result in this section is the following proposition.
\bp{prop-encage}
Fix a circuit ${\bf z}\in \E({\bf k})$ with ${\bf k}\in V(\Lambda,n)$.
For any $\epsilon>0$, there is $L$ integer, and a constant $\beta>0$
independent of $\epsilon$, such that if $\C_i:=\C(z(i))$, and
\be{def-encage}
\begin{split}
P^{L}_{z(i)}\pare{S_T=z(i+1)}=&
\ind_{\acc{z(i+1)\in \C_i}}P_{z(i)}\pare{S_T=z(i+1),T<T(\S(\C_i))}\\
&\quad+ \ind_{\acc{z(i+1)\not\in \C_i}}P_{z(i)}\pare{S_T=z(i+1)},
\end{split}
\ee
then
\be{cor-cage}
\prod_{i=0}^{|{\bf k}|-1} P_{z(i)}\pare{S_T=z(i+1)}
\le e^{\beta \epsilon {\sqrt n}} \prod_{i=0}^{|{\bf k}|-1} 
P^{L}_{z(i)}\pare{S_T=z(i+1)}.
\ee
\ep
\br{rem-encage}
Consider a $\C$-loop, say $\L$, and assume that for some integer
$i$, $\L$ corresponds to the $i$-th $\C$-loop in circuit ${\bf z}$.
We use the shorthand notation $\Weight(\L)$ 
to denote the probability associated with $\L$ 
\be{not-loop}
\Weight(\L):= \prod_{k=\tau_i}^{\sigma_i} P^L_{z(k-1)}\pare{S_T=z(k)}.
\ee
Note that $\Weight(\L)$ includes the probabilities of the
entering and exiting {\it trip}. The point of encaging loop is
the following identity
\[
\prod_{k=\tau_i}^{\sigma_i-2} P^L_{z(k)}\pare{S_T=z(k+1)}=
\prod_{k=\tau_i}^{\sigma_i-2} P^L_{\T(z(k))}\pare{S_T=\T(z(k+1))}.
\]
Thus, if we set $z=z(\tau_i-1)$ and $z'=z(\sigma_i)$
\be{encage-key}
\Weight(\L):=\frac{P_{z}\pare{S_T=z(\tau_i)}}
{P_{z}\pare{S_T=\T(z(\tau_i))}}\frac{P_{z(\sigma_i-1)}\pare{S_T=z'}}
{P_{\T(z(\sigma_i-1))}\pare{S_T=z'}}\ 
\Weight(\T(\L)).
\ee
\er
The proof of Proposition~\ref{prop-encage} is divided in two lemmas.
The first lemma deals with excursions between {\it close} sites. 
Such excursions are abundant. The larger $L$ is,
the better the estimate \reff{encage.1} of Lemma~\ref{lem-encageS}.
The second result, Lemma \ref{lem-encageB},
deals with excursions between {\it distant} sites of the
same cluster. Such excursions are rare, and even a large constant in
the bound \reff{encage.2} is innocuous.
\bl{lem-encageS}
For any $\epsilon>0$, there is $L$, such that for any $L$-cluster $\C$, and $x,y\in \C$, with
$|x-y|\le \sqrt{L}$, we have
\be{encage.1}
P_x(S_T=y)\le e^{\e} P_x\pare{ S_T=y,T<T(\S)}.
\ee
\el
\bl{lem-encageB}
There is $C_B$ independent of $L$, such that for any $L$-cluster $\C$, and $x,y\in \C$, with
$|x-y|> \sqrt{L}$, we have
\be{encage.2}
P_x(S_T=y)\le C_B P_x\pare{ S_T=y,T<T(\S)}.
\ee
\el
Lemmas~\ref{lem-encageS} and~\ref{lem-encageB} are proved in the
Appendix. We explain how they yield \reff{cor-cage}, that is
how to bound the cost of encaging a {\it loop}. 
Consider a circuit associated with ${\bf k}\in V(\Lambda',n)$ and
${\bf z}\in \E({\bf k})$.
\begin{itemize}
\item[(i)] Each journey between
sites at a distance less than ${\sqrt L}$ brings a cost $e^{\epsilon}$ from \reff{encage.1},
and even if ${\bf z}$ consisted only of such journeys, the cost would be negligible, 
since the total number of visits of $\Lambda$ is $|{\bf k}|\le 2A^4 {\sqrt n}$ as seen
in \reff{main-estimate}.
\item[(ii)] Each journey between sites at a distance 
larger than ${\sqrt L}$ brings a constant $C_B$,
but their total number is less than $\epsilon{\sqrt n}$ by the second
constrain in \reff{circuit.6}. 
\end{itemize}
Combining (i) and (ii), we obtain \reff{cor-cage}.

\subsection{Local Circuits Surgery.}\label{sec-movetrip}
In this section, we first estimate the cost of wiring differently
trips. More precisely, we have the following two lemmas.

\bl{trans-lem.1}There is a constant $C_T>0$, such that for any 
$y\in \Lambda\bs \C$ and $x\in \C$, we have
\be{trans.1}
P_y( S_T=x)\le C_T P_{y}( S_T=\T(x)).
\ee
\el
\br{rem-trans.1} By noting that for any $x,y\in \Lambda$, $P_x(S_T=y)=
P_y(S_T=x)$, we have also \reff{trans.1} with the r\^ole of $x$ and $y$ 
interchanged. However, it is important to see that the following 
inequality with $C$ independent of $n$
\be{wrong}
P_y( S_T=\T(x))\le C P_{y}( S_T=x)\quad\text{ is wrong !}
\ee
Indeed, the distance between $y$ and $\T(x)$ 
might be considerably shorter than
the distance between $y$ and $x$, 
and the constant $C$ in \reff{wrong}
should depend on this ratio of distances, and thus on $n$.
\er
Secondly, we need to wire different points of 
the same cluster to an outside point.
\bl{improp-lem.1} There is a constant $C_I>0$, such that for 
all $x,x'\in \C$, and for $y\in \Lambda'\bs \C$
\be{improp.21}
P_y(S_T=x)\le C_I P_y(S_T=x'),\text{ and for }y\in \Lambda'\bs \T(\C),
\ P_y(S_T=\T(x))\le C_I P_y(S_T=\T(x')).
\ee
Moreover, \reff{improp.21} holds when we interchange
initial and final conditions.
\el

Finally, we compare the cost of different
trips joining $\C$ and $\T(\C)$. This is a corollary of 
Lemma~\ref{improp-lem.1}.
\bc{improp-cor.1} For all $x,x'\in \C$ and $y,y'\in \C$,
\be{improp.5}
P_{x}(S_T=\T(y))\le C_I^2 P_{x'}(S_T=\T(y'))
,\quad\text{ and }\quad
P_{\T(x)}(S_T=y)\le C_I^2 P_{\T(x')}(S_T=y').
\ee
\ec

\section{Global Circuits Surgery.}\label{sec-moveloop}
In this section, we discuss the following key result.
We use the notations of 
Section~\ref{sec-defcir}.
\bp{synth-prop.1} 
There is $\beta>0$, such that for any $\epsilon>0$, 
\be{synth.1}
P\pare{||\ind_{\Lambda}l_\infty||_2^2\ge n\xi,\ \Lambda
\subset \D_\infty(A,\sqrt n)}\le
e^{\beta\epsilon\sqrt{n}}
P\pare{ ||\ind_{\T(\Lambda)}l_{\infty}||_2^2\ge n\xi,
\ \T(\Lambda) \subset \D_\infty(A,\sqrt n)}.
\ee
\ep
We iterate a finite number of times 
Proposition~\ref{synth-prop.1}, with starting set $\T(\Lambda)$,
then $\T^2(\Lambda)$ and so forth
(at most $|\Lambda|$-iterations are enough), and end up with a
finite set $\tilde \Lambda$ made up of just one $L$-cluster.

If $\dist(0,\tilde \Lambda)$ is larger than $2\diam(\tilde \Lambda)$,
then we can choose an arbitrary point $z^*$ at a 
distance $\diam(\tilde \Lambda)$ 
from $\tilde \Lambda$, and replace in the circuit
decomposition of \reff{circuit.14} 
$P_0(S_T=z(1))$, for any $z(1)\in \tilde \Lambda$,
by $P_{z^*}(S_T=z(1))$ 
at the cost of a constant, by arguments similar to those
of Section~\ref{sec-movetrip}, 
and then use translation invariance to translate
$\tilde \Lambda$ by $z^*$ back to the origin.
Thus, from Proposition~\ref{synth-prop.1}, we obtain easily the following 
result.
\bp{prop-synth2}
There is $\tilde \Lambda\ni 0$ a subset of $\Z^d$ 
whose diameter depends on $\e$ but not on $n$, 
such that for $n$ large enough
\be{synth.14}
P\pare{||\ind_{\Lambda}l_\infty||_2^2\ge n\xi,\ \Lambda
\subset \D_\infty(A,\sqrt n)}\le
e^{\beta\epsilon \sqrt{n}}
P_0\pare{||\ind_{\tilde \Lambda} l_{\infty}||_2^2\ge n\xi,
\ \tilde \Lambda \subset \D_\infty(A,\sqrt n)}.
\ee
\ep
\noindent{\bf First steps of proof of Proposition~\ref{synth-prop.1}}
Fix $\epsilon>0$. Proposition~\ref{prop-encage} produces a scale
$L$ which defines $L$-clusters, which in turn allows us to define
{\it circuits}. Also, the constant $\beta$ in \reff{cor-cage} is
independent of $\epsilon$. Recalling
\reff{circuit.14} together with \reff{cor-cage}, we obtain
\be{circuit-cage}
P\pare{||\ind_{\Lambda}l_\infty||_2^2\ge n\xi,\ \Lambda
\subset \D_\infty(A,\sqrt n)}\le
e^{\beta\epsilon \sqrt n} \sum_{{\bf k}\in V(\Lambda',n)}
\sum_{{\bf z}\in\E({\bf k})} \prod_{i=1}^{|{\bf k}|}
P^L_{z(i-1)}\pare{S_T=z(i)}.
\ee

Recall that for ${\bf k}\in V(\Lambda',n)$, $\E({\bf k})$ 
is the collection of
possible circuits producing local times ${\bf k}$ with
$\{\sum_{\Lambda} k(x)^2\ge n\xi\}$. The aim of this section
is to modifiy the circuits so as to interchange the r\^ole of
$\C_1$ and $\T(\C_1)$.

We aim at building a map $f$
on circuits with the following three properties: if ${\bf z}
\in \E({\bf k})$
\be{feat-0}
\text{(i)}\qquad\forall x\in \Lambda\bs\C,
\ l_{f({\bf z})}(x)=k(x),\ \forall x\in \C_1,\
l_{f({\bf z})}(\T(x))\ge k(x), \text{ and }
\ l_{f({\bf z})}(x)\le k(\T(x)).
\ee
Secondly, for $\beta>0$ and a constant $C(\Lambda)>0$ depeding only 
on $|\Lambda|$,
\be{feat-1}
\text{(ii)}\qquad \forall z\in f(\E({\bf k})),\quad 
|f^{-1}(z)|\le C(\Lambda) e^{\beta \epsilon \sqrt n},
\ee
Thirdly,
\be{feat-2}
\text{(iii)}\qquad
\prod_{i=0}^{|{\bf k}|-1}
P^{L}_{z(i)}\pare{S_T=z(i+1)}
\le e^{\beta \epsilon {\sqrt n}} \prod_{i=0}^{|{\bf k}|-1}
P^{L}_{f(z(i))}\pare{S_T=f(z(i+1))}.
\ee
Assume, for a moment, that we have $f$ with (i),(ii) and (iii).
Then, summing over ${\bf z}\in\E({\bf k})$, 
\be{synth.12}
\begin{split}
\sum_{{\bf z}\in\E({\bf k})}&\prod_{i=0}^{|{\bf k}|-1}
P^{L}_{z(i)}\pare{S_T=z(i+1)}
\le e^{\beta\epsilon\sqrt{n}}\sum_{{\bf z}\in\E({\bf k})}
 \prod_{i=0}^{|{\bf k}|-1}
P^{L}_{f(z(i))}\pare{S_T=f(z(i+1))}\\
&\le  e^{\beta\epsilon\sqrt{n}}\sum_{{\bf z}\in  f(\E({\bf k}))}
|f^{-1}(z)|\prod_{i=0}^{|{\bf k}|-1} 
P^{L}_{z(i)}\pare{S_T=z(i+1)}\\
&\le C(\Lambda) e^{2\beta\epsilon\sqrt{n}}\ 
P_0\pare{l_{\infty}|_{\Lambda\bs\C}=k|_{\Lambda\bs\C}
,\ \forall x\in \C_1,\ l_{\infty}(\T(x))\ge k(x),\text{ and }
\ l_{\infty}(x)\le k(\T(x))}
\end{split}
\ee
We further sum over 
${\bf k}\in V(\Lambda',n)$, and replace the sum over the
$\acc{k(y)\le A\sqrt n,\ y\in \T(\C_1)}$ 
by a factor $(A\sqrt{n})^{|\Lambda|}$,
and rearrange the sum over $\acc{k(y) ,\ y\in \C_1}$, to obtain
\be{synth.13}
\begin{split}
\sum_{{\bf k}\in V(\Lambda',n)}\sum_{{\bf z}\in \E({\bf k})} 
\prod_{i=0}^{|{\bf k}|-1}&P^{L}_{z(i)}\pare{S_T=z(i+1)}
\le e^{2\beta\epsilon\sqrt{n}} (A\sqrt{n})^{|\Lambda|}\\
&\times E\cro{\prod_{y\in \C_1} l_\infty(\T(y)),\ \T(\Lambda)
\subset \D_\infty(A,\sqrt n),
||\ind_{\T(\Lambda)}l_{\infty}||_2^2\ge n\xi}.
\end{split}
\ee
Note that in \reff{synth.13}, we can
assume $l_\infty(\T(y))\le A\sqrt{n}$
for all $y\in \C$, since for a transient walk,
the number of visits to a given site is bounded by a 
geometric random variable. Thus, in the
expectation of \reff{synth.13}, we bound $l_\infty(\tilde y)$ 
by $A\sqrt{n}$, and $|\C|$ by $|\Lambda|$. 

Providing we can show the existence of a map $f$ with
properties \reff{feat-0}, \reff{feat-1} and \reff{feat-2}, 
we would have proved Proposition \ref{synth-prop.1}.
Sections~\ref{sec-marriage}, \ref{sec-proper} and \ref{sec-improper}
are devoted to contructing the map $f$.

\subsection{A Marriage Theorem.}\label{sec-marriage}
This section deals with {\it global} modifications of circuits.
For this purpose, we rely on an old 
{\it Marriage Theorem} (see e.g.\cite{JUKNA}), which seems to have
been first proved by Frobenius~\cite{FROBENIUS} in our setting. 
Since we rely heavily
on this classical result, we quote it for the ease of reading.
\bt{circuit-th.1}{Frobenius' Theorem.} Let $\G=(G,E)$ 
be a k-regular bipartite graph
with bipartition $G_1,G_2$. Then, there is a bijection $\v:G_1\to G_2$ such that
$\acc{(x,\v(x)),\ x\in G_1}\subset E$.
\et
Now, to see how we use Frobenius' Theorem, we need more notations.
First, for two integers $n$ and $m$, we call
\be{circuit.11}
\Omega_{n,m}=\acc{\eta\in \acc{0,1}^{n+m}:\ \sum_{i=1}^{n+m} \eta(i)=n}.
\ee
Now, when $n>m$, we define the graph $\G_{n,m}=
(G_{n,m},E_{n,m})$ with
$G_{n,m}=\Omega_{n,m}\cup \Omega_{m,n}$, and
\be{circuit.12}
E_{n,m}=\acc{(\eta,\zeta)\in \Omega_{n,m}\times \Omega_{m,n}
:\ \zeta(x)\le \eta(x),\forall x\le n+m}.
\ee
With $k=n-m$, $\G_{n,m}$ is a $k$-regular graph with bipartition 
$\Omega_{n,m},\Omega_{m,n}$, and
Frobenius' Theorem gives us a bijection $\v_{n,m}:
\Omega_{n,m}\to \Omega_{m,n}$. Thus, under the action of
$\v_{n,m}$ a 1 can become a 0, but a 0 stays 0. The importance
of this feature is explained below in Remark~\ref{rem-frob}.
When $n=m$, we call $\v_{n,n}$ the identity on $\Omega_{n,n}$.

We use Frobenius' Theorem to select pairs of {\it trips} with
the same {\it type}, one {\it trip} to $\C$ 
and one {\it trip} to $\tilde \C$ which are interchanged. 
Then, we describe how the associated loops are interchanged.
However, some {\it patterns} of loops cannot be handled using
Frobenius' Theorem, and we call these loops {\it improper}.
For the ease of notations, we call $\C=\C_1$ and $\tilde \C=\T(\C_1)$.
\bd{def-proper}
A $\C$-loop is called {\it proper} if it is preceded by a trip from
$\Lambda$ to $\C$, and the other $\C$-loops are called {\it improper}.
Similarly, a $\tilde\C$-loop is called {\it proper} if it is 
preceded by a trip from $\T(\Lambda)$ to $\tilde\C$.
\ed

We describe in the two next sections, how to define a map
$f$ satisfying \reff{feat-0},\reff{feat-1} and \reff{feat-2}.
This map only transforms $\C$ and $\tilde \C$-loops.
It acts on each {\it proper} loop of
a certain {\it type}, say $p$ and $\tilde p$, by a global action that we
denote $f_p$. Also, there will be an action $f_i$ on 
{\it improper }loops which we describe in Section~\ref{sec-improper}.
Thus, $f$ is a composition of $\{f_p,p\in \C^2\}$ and $f_i$, taken in the
the order we wish. Note that for any ${\bf z}\in f(\E({\bf k}))$, we have
\be{def-f}
|f^{-1}({\bf z})|=
\prod_{p\in \C^2}|f_p^{-1}({\bf z})|\times |f_i^{-1}({\bf z})|.
\ee
Thus, property \reff{feat-1} holds for $f$, if it holds for $f_i$, and
for each $f_p$ as $p\in \C^2$.
We describe the $\{f_p,p\in \C^2\}$ in Section~\ref{sec-proper},
and $f_i$ in Section~\ref{sec-improper}.
\subsection{{\it Proper} Loops.}\label{sec-proper}
We fix ${\bf k}\in V(\Lambda',n)$ and ${\bf z}\in \E({\bf k})$.
We fix a {\it type} $p=(z,z')\in \C^2$, and
we call $\nu(p)$ the number of proper $\C$-loops of type $p$ in ${\bf z}$.
Similarly, $\nu(\tilde p)$ is the number  of proper $\tilde \C$-loops 
of type $\tilde p=(\T(z),\T(z'))$.
To each {\it type} $p$ corresponds a
configuration $\eta_p\in \Omega_{\nu(p),\nu(\tilde p)}$ which encodes
the successive occurrences of proper $\C$ and $\tilde \C$-loops 
of type $p$:
a mark 1 for a $ \C$-loop and a mark 0 for a $\tilde \C$-loop.

Assume that $n:=\nu(p)\ge m:=\nu(\tilde p)$, and $\eta_p\in \Omega_{n,m}$.
All $\C$-loop ({\it proper} and of {\it type} $p$) are translated
by $\T$, and all $\tilde \C$-loop
({\it proper} and of {\it type} $\tilde p$) are translated by $\T^{-1}$.
The bijection $\v_{n,m}$ encodes the positions of 
the translated loops, as follows.
\begin{itemize}
\item The $\C$-loop associated with the $i$-th
occurrence of a 1 in $\eta_p$, is transformed into a $\tilde\C$-loop
associated with the $i$-th occurrence of a 0 in $\v_{n,m}(\eta_p)$.
\item The $\tilde \C$-loop associated with the $i$-th
occurrence of a 0 in $\eta_p$, is transformed into a $\C$-loop
associated with the $i$-th occurrence of a 1 in  $\v_{n,m}(\eta_p)$.
\end{itemize}
After acting with $f_p$, 
the number of $\tilde\C$-loops of {\it type }$\tilde p$ 
increases by $\nu(p)-\nu(\tilde p)\ge 0$.

\begin{figure}[htpb]
\centering
\includegraphics[width=10cm,height=10cm]{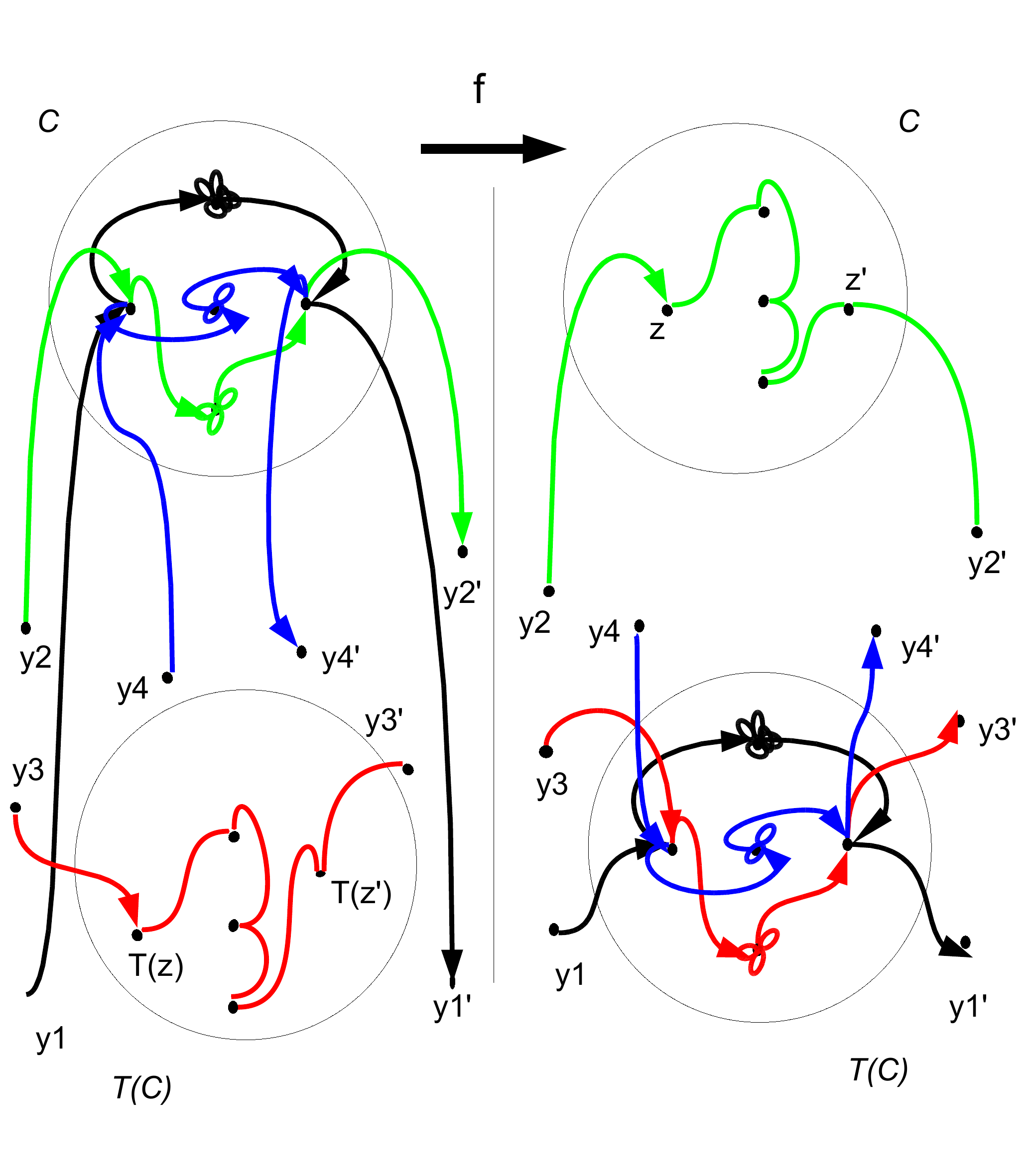}
\caption{Action of $f$ on proper loops.}\label{fig:dessin1}
\end{figure}

For definiteness, we illustrate this algorithm on a simple
example (see Figure~\ref{fig:dessin1}. 
Assume that circuit ${\bf z}\in \E({\bf k})$
has 3 proper $\C$-loops of type $p$, say $\L_1,\L_2$ and $\L_3$,
and 1 proper $\tilde \C$-loop of type $p$, say $\tilde \L_1$. Let us make visible 
in ${\bf z}$ only these very loops and the trips joining them:
\be{ex-1}
{\bf z}:\qquad\dots y_1\L_1y_1'\dots y_2\L_2y_2'\dots y_3\tilde \L_1 y_3'\dots 
y_4\L_3y_4'\dots ,
\ee
for $\{y_i,y_i',\ i=1,\dots,4\}$ in $\Lambda\bs\C$.
For such a circuit, we would have $\nu(p)=3$ and $\nu(\tilde p)=1$ and
$\eta_p=(1101)$. Furthermore, assume that $\v_{3}(1101)=0100$. 
Then, the $p,\tilde p$ proper loops are transformed into 
\be{ex-2}
{\bf f_p(z)}:
\qquad\dots y_1\T(\L_1)y_1'\dots y_2\T^{-1}(\tilde \L_1)y_2'\dots 
y_3\T(\L_2) y_3' \dots y_4\T(\L_3)y_4'\dots 
\ee
We end up with 3 $\tilde \C$-loops of type $\tilde p$,
$\T(\L_1),\T(\L_2)$ and $\T(\L_3)$, 
and one $\C$-loop $\T^{-1}(\tilde \L_1)$.
Note that in both ${\bf z}$ and ${\bf f_p(z)}$, the second loop (of type
$p$ or $\tilde p$) is a $\C$-loop, as required by Frobenius map $\v_{3}$.
The configuration $z$ in \reff{ex-1} is represented on the
left hand side of Figure~\ref{fig:dessin1}, whereas $f_p(z)$ is shown 
on its right hand side. Note that we put most
of the sites $\{y_i,y'_i,i=1,\dots,4\}$ close to $\T(\C)$. This
is the desired feature of $\T$ as established in Lemma~\ref{lem-move}.

\br{rem-frob}
One implication of the key feature of $\v_{n,m}$,
namely that $(\eta_p,\v_{n,m}(\eta_p))
\in E_{n,m}$, is that a trip $(y,\T(z))$ or $(\T(z'),y')$ is invariant
under $f_p$. Note that in Figure~\ref{fig:dessin1}, $(y_3,\T(z))$
and $(\T(z'),y_3')$ are invariant, whereas $(y_1,z)$ becomes
$(y_1,\T(z))$ and fortunately $|y_1-\T(z)|\le |y_1-z|$ on the drawing.
\er
Note that $f_p$ satisfies \reff{feat-0}. Indeed,
if we call $z_p$ the substring of $z$ made up of only sites
represented in \reff{ex-1}, and $f_p(z_p)$ the
substring of $f_p(z)$ made up of only sites
represented in \reff{ex-2}, we have $l_{f_p(z_p)}(x)=l_{z_p}(x)$ for
$x\in \Lambda\bs\C$,
\be{primer-1}
\forall x\in \C,\quad l_{f_p(z_p)}(\T(x))=l_{z_p}(x),\quad\text{and}
\quad l_{f_p(z_p)}(x)=l_{z_p}(\T(x)).
\ee
%About cost.
Now, we estimate the cost of going from $z_p$ to $f_p(z_p)$.
We consider encaged loops as described in Section~\ref{sec-cageloop}.
The purpose of having defined {\it types}, and of having
{\it encaged} loops, is the following two simple observations, 
which we deduce from \reff{encage-key} in Remark~\ref{rem-encage}.
\be{key1-type}
\text{(i)}\qquad
\Weight(\tilde \L_1)\Weight(\L_2)=
\Weight(\T^{-1}(\tilde \L_1))\Weight(\T(\L_2)),
\ee
and, if $p=(z,z')\in \C^2$
\be{loop-exchange}
\text{(ii)}\qquad
\Weight(\L_1)
=\frac{P_{y_1}\pare{S_T=z}}{P_{y_1}\pare{S_T=\T(z)}}
\frac{P_{z'}\pare{ S_T=y'_1}}{P_{\T(z')}
\pare{ S_T=y'_1}}\quad \Weight(\T(\L_1)),
\ee
and a similar equality linking $\Weight(\L_3)$ and $\Weight(\T(\L_3))$.
Thus, the cost of transformation \reff{ex-2} is $C_T^4$, 
where $C_T$ appears in 
Lemma~\ref{trans-lem.1}, since only 2 entering trips and 2 exiting
trips have been wired differently.

Now, for any ${\bf z}\in \E({\bf k})$, 
the number of loops which undergo a transformation
is less than the total number of loops, 
which is bounded by $\epsilon {\sqrt n}$.
The maximum cost (maximum over ${\bf z}\in \E({\bf k})$)
of such an operation is $2C_T$ to the power $\epsilon {\sqrt n}$.

The case (rare but possible) where $\nu(p)< \nu(\tilde p)$ 
has to be dealt with differently. Indeed,
for an arbitrary cluster $\C'$, we cannot transform a trip between 
$\C'$ and $\tilde \C$ into a trip between $\C'$ and $\C$
at a constant cost, since $\dist(\C',\tilde \C)$ 
might be much smaller than $\dist(\C',\C)$.

We propose that $f_p$ performs the following changes:
\begin{itemize}
\item Act with $\T$ on all $\C$-loops of {\it type} $p$.
\item Act with $\T^{-1}$ only on the first $\nu(p)$ $\tilde \C$-loops 
of {\it type} $\tilde p$.
\item Interchange the position of 
the $\nu(p)$ first $\C$-loops with $\nu(p)$ first $\tilde \C$-loops.
\end{itemize}
For instance, in the following example, ${\bf z}$ has
three $\tilde \C$-loops $\tilde \L_1,\tilde \L_2$ and $\tilde\L_3$ and
one $\C$-loop $\L_1$, 
\[
{\bf z}:\dots y_1\tilde \L_1y_1'\dots y_2\L_1y_2'\dots 
y_3\tilde\L_2 y_3'\dots y_4\tilde\L_3y_4'.
\]
$\nu(p)=1<\nu(\tilde p)=3$, and we have
\be{ex-3}
{\bf z}: \longrightarrow
f_p({\bf z}):\dots y_1\T(\L_1)y_1'\dots 
y_2\T^{-1}(\tilde \L_1)y_2'\dots y_3\tilde\L_2 y_3'\dots
y_4\tilde\L_3y_4'.
\ee
In so doing, note that the cost is 1, but
instead of \reff{primer-1}, we have 
\be{primer-2}
\forall x\in \C,\quad l_{f_p(z_p)}(\T(x))\ge l_{z_p}(x),\quad\text{and}
\quad \forall x\in \C,\quad l_{f_p(z_p)}(x)\le l_{z_p}(\T(x)).
\ee
Also, we have brought a multiplicity of pre-images.
Indeed, note that the final circuit of \reff{ex-3} 
could have been obtained, following
the rule of \reff{ex-2}, 
by a circuit ${\bf z'}$ where $\nu(p)\ge \nu(\tilde p)$:
\be{ex-4}
{\bf z'}:\dots y_1\L_1y_1'\dots
y_2\T^{-1}(\L_2')y_2'\dots y_3\tilde \L_1 y_3'\dots
y_4\T^{-1}(\L_3')y_4'\dots \longrightarrow f_p({\bf z}).
\ee
Also, $f_p$ maps a {\it proper} loop into a {\it proper} loop, 
and a pre-image under $f_p$ has either $\nu(p)\ge \nu(\tilde p)$ or 
$\nu(p)< \nu(\tilde p)$, and so only two possible pre-images.
Since this is true for any {\it type}, 
an upper bound on the number of pre-images of the composition of all $f_p$,
is bounded by 2 to the power $|\C|^2$ (which is the number of {\it types}).
Since $\C\subset \Lambda$ whose volume
is independent of $n$, the multiplicity is innocuous in this case.
\subsection{{\it Improper} Loops.}\label{sec-improper}
In this section, we deal with 
trips in $\C\times \tilde\C\cup \tilde \C\times \C$. The notion
of {\it type} is not useful here. We call $f_i$ the action of $f$
on {\it improper} loops.

To grasp the need to distinguish {\it proper} loops from
{\it improper} loops, assume that we have a trip from a $\T(\C)$-loop 
to a $\C$-loop. If we could allow the $\C$-loop 
to become a $\T(\C)$-loop, we 
could reach a situation with two successive $\T(\C)$-loops
linked with no trip. They would merge
into one $\T(\C)$-loop by our definition \ref{def-trip}. 
This may increase dramatically the number of pre-images of a given
$f({\bf z})$, violating \reff{feat-1}. 
\begin{figure}[htpb]
\centering
\includegraphics[width=10cm,height=10cm]{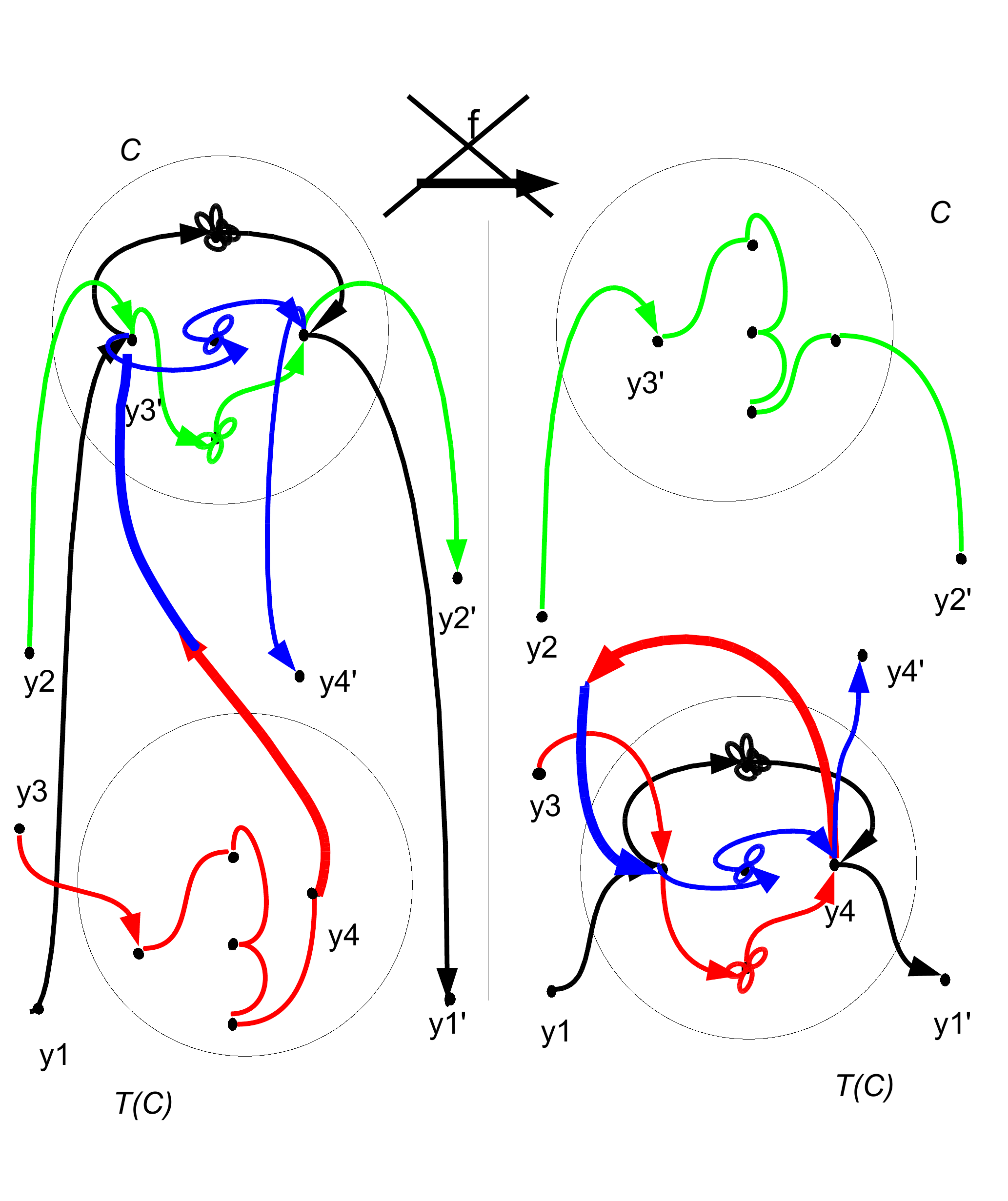}
\caption{Red and blue loops merging.}\label{fig:dessin2}
\end{figure}
We illustrate this with a concrete example drawn in 
Figure~\ref{fig:dessin2}, below.
We have considered the same example as in \reff{ex-1}, but now 
there is a trip between $\tilde \L_1$ to $\L_3$, so that
$y'_3$ is in loop $\L_3$ whereas $y_4\in \tilde \L_1$,
as shown in Figure~\ref{fig:dessin2}. If 
we where to apply the algorithm of Section~\ref{sec-proper},
we would obtain the image shown on the right hand side
of in Figure~\ref{fig:dessin2}. There, the loops $\T(\L_2)$ and
$\T(\L_3)$ (that we obtain in \reff{ex-2}) would have to merge.

Consider first
a circuit with a string of successive {\it improper loops} of type $p$, 
such that the number of $\C$-loops matches the number of $\tilde\C$-loops. 
For instance, 
assume that the $i$-th $\tilde \C$-loop 
is {\it improper} and followed by the $j$-th $\C$-loop, and so
forth. For definiteness, assume that ${\bf z}$ contains $z_i$ ($i$ for
{\it improper}) with
\be{improp.1}
z_i:=y_1\tilde \L(i)\L(j)\dots \tilde\L(i+k) \L(j+k) 
y'_1,\quad\text{with}\ k\ge 0,
\ \text{ and}\quad y_1,y'_1\not\in \C\cup\tilde \C.
\ee
Our purpose is to transform such a sequence of 
alternating $\C$-$\tilde\C$ loops into
a similar alternating sequence, such that $f_i(z_i)$ 
satisfies \reff{feat-0}, \reff{feat-1} and \reff{feat-2}.

One constraint is that we cannot replace the entering trip, and exiting trip
in general, which in turn fixes the order of visits to $\C$ and $\tilde \C$.
Indeed, as in the previous section, if $p=(z,z')$ and
$|y_1-\T(z)|\ll |y_1-z|$, then we cannot map the trip $(y_1,\T(z))$
to $(y_1,z)$ at a small cost. We propose to following map
\be{improp.2}
f_i(z_i):=y_1\T(\L(j))\T^{-1}(\tilde \L(i))\dots \T(\L(j+k))
\T^{-1}(\L(i+k)) y'_1
\ee
Note that \reff{primer-1} holds.
With an abuse of notations we represent the probability associated 
with $f_i(z_i)$, as
\be{not-improper}
\Weight(f_i(z_i)):=\prod_{l=0}^{k-1}\Weight(\tilde \L(i+l))\Weight(\L(j+l)),
\ee
even though we mean now that the trips joining successives journeys between
$\C$-$\tilde\C$ or $\tilde\C$-$\C$are counted only once.
Thus, the estimates we need concern {\it trips} 
joining {\it improper loops} together,
in addition to the first entering and the last exiting trip from $\gamma$. 
These estimates are the content of Lemma~\ref{improp-lem.1}. 
The cost $\Weight(z_i)/\Weight(f_i(z_i))$ is bounded by
$C_I^{2(k+1)+1}$, where $k+1$ is the number of successive blocks of $\tilde \C$-$\C$ loops.
Since the total number of {improper} loops of all {\it types}
is bounded by $\epsilon{\sqrt n}$, the total cost
is negligible in our order of asymptotics.

The case where the number of $\C$ and $\tilde \C$-loops
does not match is trickier. First, assume that we deal with
\be{ex-10}
z_i:=y_1\L(i)\tilde \L(j)\dots\L(i+k) y'_1.
\ee
Here, we have no choice but to replace $z_i$ with
\be{ex-11}
f_i(z_i):=y_1\T(\L(i))\T^{-1}(\L(j))\dots\T(\L(i+k)) y'_1.
\ee
Note that \reff{primer-2} holds. 

Lastly, consider the case with more $\tilde\C$-loops. For instance,
\be{ex-12}
z_i:=y_1\tilde\L(i)\L(j)\tilde\L(i+1)\dots\L(j+k-1)\tilde \L(i+k) y'_1.
\ee
For reasons already mentioned, 
we cannot map the first $\tilde \C$-loop into a $\C$-loop.
We propose to keep the first loop unchanged, 
and act on the remaining loops, in the following way
\be{ex-13}
f_i(z_i):=y_1\tilde\L(i)\T^{-1}(\tilde \L_{i+1}) \T(\L(j))
\dots \T^{-1}(\tilde \L_{i+k}) \T(\L(j+k-1))y'_1.
\ee
Here, as in \reff{ex-3}, \reff{primer-2} holds, and this choice brings
a multiplicity of pre-images. Indeed, $f_i(z_i)$ could have come from 
\[
z_i':=y_1\T^{-1}(\tilde \L(i))\tilde \L(i+1)\L(j) \dots
\tilde \L_{i+k} \T(\L(j+k-1) y'_1\longrightarrow f_i(z_i). 
\]
So, in estimating the number of pre-images
of a circuit, 
we find that it is at most 2 to the power of the number of 
{\it improper} loops.
Now, the maximum number of {\it improper} loops is $\epsilon{\sqrt n}$. 
Also, the cost of transforming all {\it improper} loops is uniformly
bounded by $C_I$ to the power $\epsilon{\sqrt n}$.

\section{Renormalizing Time.}\label{sec-time}

In this section, we show the following result.
\bp{prop-time.1} For any finite domain $\tilde \Lambda\subset\Z^d$, 
there are positive constants $\alpha_0$, $\gamma$, 
such that for any large integer $n$, there is a sequence 
${\bf k_n^*}=\acc{k_n^*(z),\ z\in \tilde \Lambda}$ with 
\be{seq-time}
\sum_{z\in \tilde \Lambda} k_n^*(z)\le n,\quad
\acc{k_n^*(z)\in [\frac{\sqrt{n}}{A}, A \sqrt{n}],
\ z\in \tilde\Lambda},\quad \text{and}\quad
\sum_{z\in \tilde\Lambda} k_n^*(z)^2\ge n\xi,
\ee
such that for any $\alpha>\alpha_0$
\be{time-main}
P_0\pare{||\ind_{\tilde \Lambda} l_{\infty}||_2^2\ge n\xi,\quad 
\tilde \Lambda\subset \D_{\infty}(A,\sqrt n)}
\le n^{\gamma} 
P_0\pare{ l_{\lfloor \alpha \sqrt{n}\rfloor }|_{\tilde\Lambda}= 
{\bf k_n^*},\ S_{\lfloor \alpha \sqrt{n}\rfloor}=0}.
\ee
\ep
\bpr
We first use a rough upper bound
\be{time.2}
P_0\pare{||\ind_{\tilde \Lambda} l_{\infty}||_2^2\ge n\xi,
\tilde \Lambda\subset \D_{\infty}(A,\sqrt n)}\le
\big|\acc{{\bf k_n}\in [\frac{\sqrt{n}}{A},
A {\sqrt n}]^{\tilde\Lambda}}\big|
\max_{{\bf k_n}\text{ in }\reff{seq-time}}
P(l_{\infty}|_{\tilde\Lambda}={\bf k_n}).
\ee
We choose a sequence ${\bf k_n^*}$
which maximizes the last term in \reff{time.2}. 
Then, we decompose $\acc{l_{\infty}|_{\tilde\Lambda}={\bf k_n^*}}$ 
into all possible circuits
in a manner similar to the circuit decomposition of 
Section~\ref{sec-circuit}:
We set $\nu=\sum_{\tilde\Lambda} k_n^*(x)$ 
(and $\nu\le |\tilde\Lambda| A\sqrt{n}$), and 
\be{time.3}
\E^*=\acc{{\bf z}=(z(1),\dots,z(\nu))
\in \tilde\Lambda^\nu:\ 
l_{{\bf z}}(x)=k_n^*(x),\ \forall x\in \tilde \Lambda}.
\ee
Then, if $T=\inf\acc{n\ge 0:\ S_n\in \tilde \Lambda}$, (and $z(0)=0$)
\be{time.4}
P_0(l_{\infty}|_{\tilde\Lambda}={\bf k_n^*})=\sum_{{\bf z}\in \E^*} 
\prod_{i=0}^{\nu-1}
P_{z(i)}\pare{\tilde T(z({i+1}))=T<\infty}P_{z(\nu)}(T=\infty).
\ee
For a fixed ${\bf z}\in \E^*$,
we call $\tau^{(i)}$ the duration of the flight 
from $z(i-1)$ and $z(i)$ which avoids
other sites of $\tilde\Lambda$. Thus,
$\tau^{(1)} \overset{\text{law} }{=} \tilde T(z(1))
\ind\{\tilde T(z(1))=T\}$, when 
restricting on the values $\acc{1,2,\dots}$, and by induction
\be{time.5}
\tau^{(i)} \overset{\text{law} }{=}  \tilde T(z(i))\circ\theta_{\tau^{(i-1)}}
\ind\acc{\tilde T(z(i))\circ\theta_{\tau^{(i-1)}}=T\circ\theta_{\tau^{(i-1)}}}
\ee
If $\TT({\bf z})=\acc{0<\tau^{(i)} <\infty,\forall i=1,\dots,\nu}$, we have
\be{time.6}
P_0(\TT({\bf z}))=
\prod_{i=0}^{\nu-1} P_{z(i)}\pare{\tilde T(z({i+1}))=T<\infty}.
\ee
Now, we fix ${\bf z}\in \E^*$ such that $P_0(\TT({\bf z}))>0$, and we fix $i<\nu$.
For ease of notations, we rename $x=z(i-1)$ and $y=z({i})$. 
Now, note that $\acc{0<\tau^{(i)}<\infty}$ contributes
to \reff{time.6} if $P_x(S_T=y)>0$, or in other words, if there is 
at least one path going from $x$ to $y$
avoiding other sites of $\tilde \Lambda$. Since $\tilde\Lambda$ has finite diameter, 
we can choose a finite length 
self-avoiding paths, and have 
\be{time.7}
P_x\pare{T(y)<\tilde T(x)}\ge 
c_\Lambda(x,y):= P_x(S_T=y,T<\infty)>c_\Lambda>0,
\ee
where $c_\Lambda$ is the minimum of $c_\Lambda(z,z')$ 
over all $z,z'\in \tilde\Lambda$ with
$c_\Lambda(z,z')>0$. Now, note that, when $S_0=y$
\be{time.8}
\tilde T(y)\ind_{T(x)<\tilde T(y)<\infty}\le \tilde T(y) \ind_{\tilde T(y)<\infty}.
\ee
Thus,
\be{time.9}
\begin{split}
E_y\cro{ \tilde T(y) \ind_{\tilde T(y)<\infty}}&\ge 
E_y\cro{\ind_{T(x)<\tilde T(y)<\infty}\pare{\tilde 
T(y)\circ\theta_{T(x)}+T(x)}}\\
&= E_y\cro{\ind_{T(x)<\tilde T(y)<\infty} T(x)}
+E_y\cro{\ind_{T(x)<\tilde T(y)<\infty}\tilde T(y)\circ\theta_{T(x)}}.
\end{split}
\ee
Now, by the strong Markov's property
\be{time.10}
E_y\cro{ \tilde T(y) \ind_{\tilde T(y)<\infty}}\ge P_y\pare{ T(x)<\tilde T(y)}
E_x\cro{T(y)\ind_{T(y)<\infty}}.
\ee
By using translation invariance of the walk and \reff{time.10}, 
we obtain
\be{time.11}
E_x\cro{ T(y)\ind_{T(y)=T<\infty}}\le E_x\cro{ T(y)\ind_{T(y)<\infty}}\le
\frac{E_0\cro{ \tilde T(0) \ind_{\tilde T(0)<\infty}}}{P_y\pare{ T(x)<\tilde T(y)}}.
\ee
Now, it is well known that there is a constant $c_d>0$ such that
for any integer $k$, $P_0(\tilde T(0) =k)\le c_d/k^{\frac{d}{2}}$,
which implies that $E_0\cro{ \tilde T(0) 
\ind_{\tilde T(0)<\infty}}<\infty$ in $d\ge 5$, and
\ba{time.12}
E_x\cro{ T(y)|T(y)=T<\infty}&=& \frac{E_x\cro{ T(y)\ind_{T(y)=T<\infty}}}{
P_x(T(y)=T<\infty)}\cr
&\le& \frac{E_0\cro{ \tilde T(0)\ind_{\tilde T(0)<\infty}}}{P_x(T(y)=T<\infty)
P_y(T(x)<\tilde T(y))}\cr
&\le & \frac{E_0\cro{ \tilde T(0) \ind_{\tilde T(0)<\infty}}}{c_\Lambda^2}.
\ea
When translating \reff{time.12} in terms of the $\acc{\tau^{(i)}}$, 
we obtain for any $\beta>0$
\be{time.13}
P\pare{ \sum_{i=1}^\nu \tau^{(i)}>\beta \nu\big| \TT({\bf z})}\le
\frac{E_0\cro{ \sum_{i=1}^\nu \tau^{(i)}|\TT({\bf z})}}{\beta \nu}\le
\frac{ E_0\cro{ \tilde T(0) \ind_{\tilde T(0)<\infty}}}{c_\Lambda^2}\times
\frac{1}{\beta}.
\ee
Thus, we can choose $\beta_0$ large enough 
(independent of ${\bf z}$) so that 
\be{time.14}
P_0\pare{ \sum_{i=1}^\nu \tau^{(i)}>\beta_0 \nu|\TT({\bf z})}\le \frac{1}{2}.
\ee
We use now
\[
P_0(\TT({\bf z}))=P_0\pare{ \sum_{i=1}^\nu \tau^{(i)}>\beta_0 \nu|\TT({\bf z})}P_0(\TT({\bf z}))+
P_0\pare{ \sum_{i=1}^\nu \tau^{(i)}\le \beta_0 \nu|\TT({\bf z})}P_0(\TT({\bf z})),
\]
to conclude that
\be{time.15}
P_0(\TT({\bf z}))\le 2 
P_0\pare{\{\sum_{i=1}^\nu \tau^{(i)}\le \beta_0 \nu\}\cap \TT({\bf z})}.
\ee
Now, there is $\alpha_0$ such that $\beta_0 \nu\le \alpha_0 \sqrt{n}$. Also, note
that there is $n_0$ such that for any $z(\nu)\in \tilde \Lambda$, there is a path
of length $n_0$ joining $z(\nu)$ to 0. Now, fix $\alpha>2\alpha_0$, 
take $n$ large enough so that 
$\lfloor \alpha \sqrt{n} \rfloor\ge \lfloor \alpha_0 
\sqrt{ n}\rfloor+ n_0$, and use
classical estimates on return probabilities,
to obtain that for a constant $C_d$ 
\be{time.30}
P_0(\TT({\bf z}))\le C_d (\alpha n)^{d/2}\!\! \sum_{\nu\le k\le \beta_0\nu}
\!\! P_0\big(\{\sum_{i=1}^\nu \tau^{(i)}=k\} \cap \TT({\bf z})\big)
P_{z(\nu)}(S_{n_0}=0)
P_0(S_{\lfloor \alpha n\rfloor -(k+n_0)}=0).
\ee
After summing over $z\in \E^*$, we obtain for any $\alpha>2\alpha_0$
\be{time.31}
\sum_{z\in \E^*}P_0(\TT({\bf z}))\le C_d (\alpha n)^{d/2} 
P_0\pare{
||\ind_{\tilde \Lambda} l_{\lfloor \alpha \sqrt{n}\rfloor }||_2^2\ge n \xi,
S_{\lfloor \alpha \sqrt{ n}\rfloor }=0}.
\ee
Note that another power of $n$ arises from the term
in \reff{time.2} yielding the desired result.
\epr
\section{Existence of a Limit.}\label{sec-exist}
We keep notations of Section~\ref{sec-time}.
We reformulate Proposition~\ref{prop-time.1} as follows.
For any finite domain $\tilde \Lambda\subset\Z^d$,
there are positive constants $\alpha_0$, $\gamma$,
such that for any $\alpha>\alpha_0$, and $n$ large
\be{start-main}
P_0\pare{||\ind_{\tilde \Lambda} l_{\infty}||_2^2\ge n\xi,\quad
\tilde \Lambda\subset \D_{\infty}(A,\sqrt n)}
\le n^{\gamma}
P_0\pare{||\ind_{\tilde \Lambda} l_{\lfloor \alpha \sqrt{n}\rfloor }
||_2\ge \sqrt{n\xi},\ S_{\lfloor \alpha \sqrt{n}\rfloor }=0}.
\ee
Thus, \reff{start-main} is the starting point in this section.
\subsection{A Subadditive Argument.}\label{sec-subadditive}
We consider a fixed region $\Lambda\ni 0$, and
first show the following lemma.
\bl{lem-sub.1}
Let $q>1$.
For any $\xi>0$ and $\Lambda$ finite subset of $\Z^d$, the 
following limit exists
\be{lawler.9}
\lim_{n\to\infty} \frac{\log(P_0(||\ind_{\Lambda}l_n||_{q}\ge n\xi,
\ S_n=0))}{n}= -I(\xi,\Lambda).
\ee
\el
\bpr 

We fix two integers $K$ and $n$, with $K$ to be taken first
to infinity. Let $m,r$ be integers such that $K=mn+r$, and $0\le r<n$.
The phenomenon behind the subadditive arguement is that
\be{lawler.1}
\A_K(\xi,\Lambda)=\acc{ ||\ind_{\Lambda}l_K||_q\ge K\xi, \ S_K=0}
\ee
is built by concatenating the ${\it same }$ optimal scenario
realizing $\A_n(\xi,\Lambda)$ on $m$ consecutive time-periods of length
$n$, and one last time-period of length $r$ where the scenario is
necessarly special and its cost innocuous. The crucial independence
between the different periods is obtained as we force the walk
to return to the origin at the end of each time period.

Our first step is to exhibit an optimal strategy
realizing $\A_n(\xi,\Lambda)$.
By optimizing over a finite number of variables
$\{k_n(x),x\in \Lambda\}$ satisfying
\be{lawler.2}
\sum_{x\in \Lambda} k_n(x)^q \ge (n\xi)^q,
\quad\text{and}\quad
\sum_{x\in \Lambda} k_n(x)\le n,
\ee
there is a 
sequence ${\bf k_n^*}:=\acc{k_n^*(x),x\in \Lambda}$ and $\gamma>0$ 
(both depend on $\Lambda$) such that
\be{lawler.3}
P_0\pare{\A_n(\xi,\Lambda)}\le n^\gamma P_0\pare{\A_n^*(\xi,\Lambda)},
,\quad\text{with}\quad
\A_n^*(\xi,\Lambda)=\acc{l_n|_{\Lambda} ={\bf k_n^*},\ S_n=0}.
\ee
Let $z^*\in \Lambda$, be the site where ${\bf k_n^*}$ reaches its 
maximum. We start witht the case $z^*=0$, and postpone the
case $z^*\not= 0$ to Remark~\ref{rem-z}.
When $z^*=0$, for any integer $r$, we call
\be{lawler.4}
\RR_r=\acc{l_r(0)=r},\quad
\text{and note that}\quad P_0(\RR_r)= P_0(S_1=0)^{r-1}>0.
\ee
Now, denote by $\A_n^{(1)},\dots,\A_n^{(m)}$
$m$ independent copies of $\A_n^*(\xi,\Lambda)$ which
we realize on the successive increments of the random walk
\[
\forall i=1,\dots,m,\quad \A_n^{(i)}=\acc{
l_{[(i-1)n,in[}|_{\Lambda}={\bf k_n^*},\ S_{in}=0}.
\]
Make a copy of $\RR_r$ independent of 
$\A_n^{(1)},\dots,\A_n^{(m)}$, by using increments
after time $nm$: that is $\RR_r=\{S_j=0,\ \forall j\in [nm,K[\}$.
Note that by independence 
\be{lawler.6}
\begin{split}
P_0\pare{\A_n(\xi,\Lambda)}^m P_0(\RR_r)\le&
(n^\gamma)^m P_0(\A_n^{(1)})\dots P_0(\A_n^{(m)})P_0(\RR_r)\\
\le &(n^\gamma)^m P_0(\bigcap_{j=1}^m \A_n^{(j)})\cap \RR_r).
\end{split}
\ee
Now, the local times is positive, so that
\[
\begin{split}
\bigcap_{i=1}^m&\acc{
l_{[(i-1)n,in[}|_{\Lambda}={\bf k_n^*},\ S_{in}=0}\cap
\acc{l_{[mn,K[}(0)=r}\\
&\subset
\acc{\sum_{x\in \Lambda}
\cro{\sum_{i=1}^m l_{[(i-1)n,in[}(x)+\ l_{[mn,K[}(x)}^q
\ge \sum_{x\in \Lambda}\pare{mk_n^*(x)+r\delta_{0}(x)}^q,S_{K}=0}.
\end{split}
\]
At this point, observe the following fact whose simple inductive
proof we omit: for $q>1$, and for $\v$ and $\psi$
are positive functions on $\Lambda$, and for $z^*\in \Lambda$,
$\v(z^*)=\max \v$, then
\be{fact-key}
(\v(z^*)+\sum_{z\in \Lambda}\psi(z))^q+\sum_{z\not= z^*}
\v(z)^q \ge \sum_{z\in\Lambda} \pare{\v(z)+\psi(z)}^q.
\ee
\reff{fact-key} implies that for any integer $m$
\ba{lawler.5}
\sum_{x\in \Lambda} \pare{ mk_n^*(x)+r\delta_{z^*}(x)}^q&\ge &
\sum_{x\in \Lambda} \pare{ mk_n^*(x)+\frac{r}{n}k_n^*(x)}^q\cr
&=& (m+\frac{r}{n})^q\sum_{x\in \Lambda} k_n^*(x)^q 
\ge (mn+r)^q \xi^q=(K\xi)^q.
\ea
Using \reff{lawler.5}, \reff{lawler.6} yields
\be{lawler.6bis}
P_0\pare{\A_n(\xi,\Lambda)}^m P_0(\RR_r)\le
(n^\gamma)^m P_0\pare{ ||\ind_{\Lambda} l_K||_q\ge K \xi,\ S_{K}=0}
\le (n^\gamma)^m P_0\pare{\A_K(\xi,\Lambda)}.
\ee
We now take the logarithm on each side of \reff{lawler.6}
\be{lawler.7}
\frac{nm}{nm+r} \frac{\log(P_0(\A_n(\xi,\Lambda)))}{n}+  
\frac{\log(P_0(\RR_r))}{K}\le 
\frac{m(\log(n^\gamma))}{nm+r}+\frac{\log(P_0(\A_K(\xi,\Lambda)))}{K}.
\ee
We take now the limit $K\to\infty$ while $n$ is kept fixed (e.g. $m\to\infty$) so that
\be{lawler.8}
\frac{\log(P_0(\A_n(\xi,\Lambda)))}{n}\le \frac{\log(n^\gamma)}{n}
+\liminf_{K\to\infty} \frac{\log(P_0(\A_K(\xi,\Lambda)))}{K}.
\ee
By taking the limit sup in \reff{lawler.8} as $n\to\infty$, we conclude that the limit 
in \reff{lawler.9} exists.

\br{rem-z}
We treat here the case $z^*\not= 0$. 
In this case, we cannot consider $\RR_r$ since
to use \reff{fact-key}, we would need the walk to start on site $z^*$,
whereas each period of length $n$ sees the walk returning to the origin.
Note that this problem
is related to the strategy on a single time-period of length $r$.
The remedy is simple: we insert a period of length $r$ into
the first time-period of length $n$ at the first time the walk hits
$z^*$; then, the walk stays at $z^*$ during $r-1$ steps.
In other words,
let $\tau^*=\inf\{n\ge 0:\ S_n=z^*\}$, $\RR_r^*=\{l_r(z^*)=r\}$ 
and note that
\be{sub.41}
\begin{split}
P_0(\A_n^{(1)})P_{z^*}(\RR_r^*)=&\sum_{i=1}^n P_0(\A_n^{(1)},\tau^*=i)
P_{z^*}(l_r(z^*)=r)\\
\le & P_0\pare{ l_{[0,n+r[}|_{\Lambda}={\bf k_n^*}+r\delta_{z^*}}.
\end{split}
\ee
Note that $P_{z^*}(\RR_r^*)=P_0(\RR_r)$, and
\[
\bigcap_{j=1}^m \A_n^{(j)}\subset
\acc{\sum_{x\in \Lambda}
\cro{l_{[0,n+r[}(x)+\sum_{i=2}^m l_{[(i-1)n,in[}(x)}^q
\ge \sum_{x\in \Lambda}\pare{mk_n^*(x)+r\delta_{z^*}(x)}^q,S_{K}=0}.
\]
We can now resume the proof of the case $z^*=0$ at step \reff{lawler.5}.
\er

\epr 
\subsection{Lower Bound in Proposition~\ref{prop-sub}.}\label{sec-LB}
We prove here the lower bound of \reff{lawler.10}.
Call $t_n$ be the integer part of $\alpha\sqrt{n}$, and
consider the following {\it scenario}
\be{lower.1}
\S_n(\Lambda,\alpha,\epsilon):=
\acc{ ||\ind_\Lambda l_{[0,t_n[}||_2^2\ge
n\xi(1+\epsilon),\ S_{t_n}=0}
\cap \acc{||l_{[t_n,n[}||_2^2
-E_0\cro{||l_n||_2^2}\ge n\xi(1-\epsilon)},
\ee
Note that $\S_n(\Lambda,\alpha,\epsilon)
\subset \{\overline{||l_n||_2^2}\ge n\xi\}$. 
Indeed, note that for any $\beta\ge 1$, and $a,b>0$ we have
$a^\beta+b^\beta\le (a+b)^\beta$. Thus, for any $x\in \Z^d$
\be{lower.3}
l^2_{[0,t_n[}(x)+ l^2_{[t_n,n[}(x)\le l^2_n(x),
\ee
and we obtain on $\S_n(\Lambda,\alpha,\epsilon)$
\be{lower.4}
E_0\cro{||l_n||_2^2}+n\xi\le 
\sum_{x\in \Lambda} l^2_{[0,t_n[}(x)+
\sum_{x\in \Z^d} l^2_{[t_n,n[}(x)
\le ||l_n||_2^2.
\ee
Note that 
$||\ind_\Lambda l_{[0,t_n[}||_2$ 
and $S_{t_n}=0$
only depend on the increments
of the random walk in the time period $[0,t_n[$, whereas
$||l_{[t_n,n[}||_2$ depends on the increments in $[t_n,n[$. Thus,
\be{lower.5}
\begin{split}
P\pare{\S_n(\Lambda,\alpha,\epsilon)}
=&P_0\pare{ ||\ind_\Lambda l_{[0,t_n[}||_2^2\ge
n\xi(1+\epsilon),\ S_{t_n}=0} \\
&\times P_0\pare{||l_{[t_n,n[}||_2^2
-E_0\cro{||l_n||_2^2}\ge n\xi(1-\epsilon)}.
\end{split}
\ee
Now, since $\frac{1}{n} ||l_n||_2^2$ converges in $L^1$ 
towards $\gamma_d$, we have $E_0[||l_n||_2^2]\le n\gamma_d(1+\epsilon/2)$
for $n$ large enough, and
we have
\be{lower.6}
P_0\pare{||l_{[t_n,n[}||_2^2
-E_0\cro{||l_n||_2^2}\ge n\xi(1-\epsilon)}\le 
P_0\pare{\frac{||l_{[0,n-t_n[}||_2^2}{n-t_n}
\ge \frac{\gamma_d-\frac{\epsilon}{2}
\xi}{1-\frac{t_n}{n}}}\longrightarrow 1.
\ee
\br{rem-LBalpha}
Note that for any $\Lambda$ finite subset of $\Z^d$, any $\beta>0$ and $\epsilon>0$
small, we have for $\chi<\zeta<1$, and $n$ large enough
\be{lower.8}
\acc{||\ind_\Lambda l_{\lfloor \beta n^\zeta\rfloor }
||_{\alpha^*}\ge \xi n^\zeta(1+\epsilon),\ 
S_{\lfloor \beta n^\zeta\rfloor}=0}
\subset \acc{||\ind_{\bar \D_n(\chi)} l_n||_{\alpha^*}\ge \xi n^\zeta}.
\ee
\er
\subsection{Proof of Theorem~\ref{intro-th.1}}\label{sec-ergodic}
First, the upper bound of Proposition~\ref{prop-sub} follows
after combining inequalities \reff{circuit.2}, \reff{circuit.14},
\reff{synth.1} and \reff{start-main}. The lower bound
of Proposition~\ref{prop-sub} is shown in the previous section.
Then, we invoke Lemma~\ref{lem-sub.1} with $q=2$, 
we take the logarithm on each sides of \reff{lawler.10},
we normalize by $\sqrt{n}$, and take the limit 
$n$ to infinity. We obtain that for any $\epsilon>0$,
there are $\alpha_\epsilon$ and $\Lambda_\epsilon$ such that 
for $\Lambda,\Lambda'\supset\Lambda_\epsilon$, and 
$\alpha,\alpha'>\alpha_\epsilon$
\be{lawler.key1}
\begin{split}
-\alpha' \ I\big(\frac{\sqrt{\xi(1+\epsilon)}}{\alpha'}&,\Lambda'\big)\le 
\liminf_{n\to\infty} 
\frac{\log\pare{ P_0(\overline{||l_n||^2_2}\ge n\xi)}}{\sqrt n}\\
&\le \limsup_{n\to\infty} 
\frac{\log\pare{ P_0(\overline{||l_n||_2^2}\ge n\xi)}}{\sqrt n}
\le -\alpha \ I\pare{\frac{\sqrt {\xi(1-\epsilon)}}{\alpha},\Lambda}+C\e.
\end{split}
\ee
By using \reff{lawler.key1}, we obtain for any $\Lambda,\Lambda'\supset 
\Lambda_\epsilon$, and $\alpha,\alpha'>\alpha_\epsilon$
\be{lawler.11}
\frac{\alpha'}{\sqrt{\xi(1+\epsilon)}} 
I\pare{\frac{\sqrt{\xi(1+\epsilon)}}{\alpha'},\Lambda'}\ge 
\sqrt{\frac{1-\epsilon}{1+\epsilon}} 
\frac{\alpha}{\sqrt{\xi(1-\epsilon)}}
I\pare{ \frac{\sqrt{\xi(1-\epsilon)}}{\alpha},\Lambda}-\frac{C\e}{\sqrt{\xi(1+\epsilon)}}.
\ee
Thus, if we call $\v(x,\Lambda)=I(x,\Lambda)/x$, 
we have: $\forall \e>0$, there is 
$x_\epsilon,\Lambda_\epsilon$ such that for $x,x'<x_\epsilon$ 
and $\Lambda,\Lambda'\supset \Lambda_\epsilon$
\be{lawler.12}
\v(x',\Lambda')\ge \sqrt{\frac{1-\epsilon}{1+\epsilon}} \quad\v(x,\Lambda)-
\frac{C\e}{\sqrt{\xi(1+\epsilon)}}.
\ee
By taking the limit $\Lambda'\nearrow\Z^d$, $x'\to 0$, 
and then $\Lambda\nearrow\Z^d$ and $x\to 0$, we reach for any $\e>0$
\be{lawler.13}
\liminf_{\Lambda\nearrow\Z^d, x\to 0} \v(x,\Lambda)\ge 
\sqrt{\frac{1-\epsilon}{1+\epsilon}} \quad\limsup_{\Lambda\nearrow\Z^d, x\to 0}  
\v(x,\Lambda)-\frac{C\e}{\sqrt{\xi(1+\epsilon)}}.
\ee
Since \reff{lawler.13} is true for $\epsilon>0$ arbitrarily small,
this implies that the limit of $\v(x,\Lambda)$ 
exists as $x$ goes to $0$ and $\Lambda$
increases toward $\Z^d$. We call this latter limit $\I(2)$, 
where the label 2
stresses that we are dealing with the $l^2$-norm of the local times.

Now, recall that the result of \cite{AC05}, (see Lemma~\ref{level-lem.1})
 says
that there are two positive constants $\underline{c},\bar c$ such that
for $x$ small enough $\underline{c}\le I(x,\Lambda)/x\le \bar c$, which together with
\reff{lawler.13} imply $0<\underline{c}\le \I(2)\le \bar c<\infty$. 
Now, using \reff{lawler.12} again, we obtain
\be{lawler.14}
\alpha  I\pare{\frac{\sqrt{\xi(1+\epsilon)}}{\alpha},\Lambda}\le \frac{1+\epsilon}
{\sqrt{1-\e}}\ \I(2)\sqrt{\xi}+C\e \sqrt{\frac{1-\epsilon}{1+\epsilon}} ,
\ee
and,
\be{lawler.15}
\alpha  I\pare{\frac{\sqrt{\xi(1-\epsilon)}}{\alpha},\Lambda}\ge 
\frac{1-\e}{\sqrt{1+\e}}\ \I(2)\sqrt{\xi}-C\e \sqrt{\frac{1-\epsilon}{1+\epsilon}}.
\ee
This establishes the Large Deviations Principle of \reff{intro.4} as $\e$ is sent to zero.
\qed

\noindent{\bf Proof of Proposition~\ref{prop-alpha}}
Looking at the proof of Theorem~\ref{intro-th.1}, we notice that the only
special feature of $\{\overline{||l_n||^2_2}\ge n\xi\}$ 
which we used, was that the excess self-intersection was realized on a 
{\bf finite} set $\D_n(A,\sqrt n)$.
Similarly, when considering 
$\{||\ind_{\bar \D_n(n^b)} l_n||_{\alpha^*}\ge \xi n^\zeta\}$,
inequality \reff{level.20} of Lemma~\ref{level-lem.3}, 
ensures that our large deviation is realized on $\D_n(A,n^\zeta)$, 
and by \reff{level.21}, we make a negligible error assuming it is not
finite. Thus, our key steps work in this case as well:
{\it circuit surgery}, {\it renormalizing time}, 
and the {\it subadditive argument}. Besides, by Remark~\ref{rem-LBalpha},
the lower bound follows trivially as well. Instead of 
\reff{lawler.10}, we would have that there is a constant $\beta$ such that
for any $\epsilon>0$, there is $\tilde \Lambda$
set of finite diameter, and $a_0>0$, such that for 
$\Lambda$ finite with $\Lambda\supset \tilde\Lambda$ and $a\ge a_0$, 
\be{lawler.30}
\begin{split}
P_0&\pare{||\ind_{\Lambda}l_{\lfloor a n^\zeta\rfloor}||_{\alpha^*}
\ge \xi(1+\e) n^\zeta, S_{\lfloor a n^\zeta\rfloor}=0 }
\le P_0\pare{||\ind_{\bar \D_n(n^b)} l_n||_{\alpha^*}\ge \xi n^\zeta}\\
&\le e^{\beta\epsilon n^\zeta} P_0\pare{||\ind_{\Lambda}
l_{\lfloor a n^\zeta\rfloor}||_{\alpha^*}\ge \xi(1-\e) n^\zeta,
S_{\lfloor a n^\zeta\rfloor}=0}.
\end{split}
\ee
Following the last step of the proof of Theorem~\ref{intro-th.1},
we prove Proposition~\ref{prop-alpha}.
\qed
\section{On Mutual Intersections.}\label{sec-KMSS}
\subsection{Proofs of Proposition~\ref{prop-CM}.}
Proposition~\ref{prop-CM} is based on the idea that
$\langle l_\infty,\tilde l_\infty\rangle$ is not {\it critical} 
in the sense
that even when {\it weighting less} intersection local times, the strategy
remains the same. In other words, define for $1<q\le 2$
\be{sinai.1}
\zeta(q)=\sum_{z\in \Z^d} l_\infty(z)\tilde l_\infty^{q-1} (z).
\ee
Then, we have the following lemma, interesting on its own.
\bl{lem-amine}
Assume that $d\ge 5$. For any $2\ge q>\frac{d}{d-2}$, there is $\kappa_q>0$ such that
\be{sinai.2}
\P\pare{\zeta(q)>t}\le \exp(-\kappa_q t^{\frac{1}{q}}).
\ee
\el
We prove Lemma~\ref{lem-amine} in the next section. 
Proposition~\ref{prop-CM} follows 
easily from Lemma~\ref{lem-amine}. Indeed, if
$\D(\xi)=\{z: \min(l_\infty(z),\tilde l_\infty(z))<\xi\}$ and $q<2$
\be{sinai.20}
\begin{split}
\acc{\bra{\ind_{\D(\epsilon \sqrt t)} l_\infty,\tilde l_\infty}>t}\subset&
\acc{\sum_{l_\infty(z)\le \frac{\sqrt t}{A}} 
l_\infty(z)^{q-1}\tilde l_\infty(z)> 
\frac{t}{2}\pare{\frac{A}{\sqrt t}}^{2-q}}\\
&\cup \acc{\sum_{\tilde l_\infty(z)\le \frac{\sqrt t}{A}} 
l_\infty(z)\tilde l_\infty(z)^{q-1}> 
\frac{t}{2}\pare{\frac{A}{\sqrt t}}^{2-q}}.
\end{split}
\ee
Then, since $1>\frac{2-q}{2}$, Lemma~\ref{lem-amine} applied to \reff{sinai.10}
implies that for large $t$
\be{sinai.11}
\P\pare{\bra{\ind_{\D(\epsilon \sqrt t)} 
l_\infty,\tilde l_\infty}>t}\le 2\exp\pare{ -\kappa_d A^{\frac{2-q}{q}}
t^{1/2}},\quad\text{since }\quad \frac{1}{q}(1-\frac{2-q}{2})=\frac{1}{2}.
\ee
\subsection{Proof of Lemma~\ref{lem-amine}.}
We assume $d\ge 5$. Lemma~\ref{lem-amine} can be thought
of as an interpolation inequality between 
Lemma 1 and Lemma 2 of \cite{KMSS}, whose proofs follow a classical pattern
(in statistical physics) of estimating all moments of $\zeta(q)$. This control is possible
since all quantities are expressed in terms of iterates of the Green's function,
whose asymptotics are well known (see for instance Theorem 1.5.4 of 
\cite{LAWLER}).

From \cite{KMSS}, it is enough that for a positive constant $C_q$,
we establish the following control on the moments
\be{sinai.3}
\forall n\in \N,\qquad \EE[\zeta(q)^n]\le  C_q^n (n!)^q.
\ee
First, noting that $q-1\le 1$, we use Jensen's inequality in the last inequality
\ba{sinai.4}
\EE[\zeta(q)^n]&\le & \sum_{z_1,\dots,z_n\in \Z^d} E_0\cro{
\prod_{i=1}^n l_\infty(z_i)}  E_0\cro{\prod_{i=1}^n l_\infty(z_i)^{q-1}}\cr
&\le & 
\sum_{z_1,\dots,z_n\in \Z^d} \pare{E_0\cro{\prod_{i=1}^n l_\infty(z_i)}}^{q}
\ea
If $\S_n$ is the set of permutation of $\acc{1,\dots,n}$ (with
the convention that for $\pi\in \S_n$, $\pi(0)=0$)
we have,
\be{sinai.5}
\begin{split}
E\cro{\prod_{i=1}^nl_\infty(z_i)}= &
\sum_{s_1,\dots,s_n\in \N} P_0(S_{s_i}=z_i,\ \forall i=1,\dots,n)\\
\le& \sum_{\pi\in \S_n} \sum_{s_1\le s_2\le\dots\le s_n\in \N}
P_0(S_{s_i}=z_{\pi(i)},\ \forall i=1,\dots,n)\\
\le & \sum_{\pi\in \S_n} \prod_{i=1}^n G_d\pare{z_{\pi(i-1)},z_{\pi(i)}} .
\end{split}
\ee
Now, by H\"older's inequality
\be{sinai.9}
\begin{split}
\sum_{z_1,\dots,z_n}\pare{\sum_{\pi\in \S_n}
\prod_{i=1}^n G_d\pare{z_{\pi(i-1)},z_{\pi(i)}}}^q&\le 
\sum_{z_1,\dots,z_n} (n!)^{q-1} \sum_{\pi\in \S_n}
\prod_{i=1}^n G_d\pare{z_{\pi(i-1)},z_{\pi(i)}}^q\\
&= (n!)^q \sum_{z_1,\dots,z_n} \prod_{i=1}^n G_d\pare{z_{i-1},z_{i}}^q.
\end{split}
\ee
Classical estimates for the Green's function, \reff{sinai.9} 
implies that
\be{sinai.10}
\begin{split}
\sum_{z_1,\dots,z_n\in \Z^d} \pare{E_0\cro{\prod_{i=1}^n l_\infty(z_i)}}^{q}\le&
(n!)^q C^n \sum_{z_1,\dots,z_n} \prod_{i=1}^n 
(1+||z_i-z_{i-1}||)^{q(2-d)}\\
\le & (n!)^q C^n \pare{ \sum_{z\in \Z^d} (1+||z||)^{q(2-d)}}^n.
\end{split} 
\ee
Thus, when $d\ge 5$ and $q>\frac{d}{d-2}$, we have a constant $C_q>0$ such
that
\be{sinai.7}
\EE[\zeta(q)^n]\le C_q^n(n!)^q
\ee
The proof concludes now by routine consideration (see e.g. \cite{KMSS} or
\cite{chen-morters}).
\subsection{Identification of the rate function \reff{identification}.}
\label{sec-ident}
The main observation is that the proof of Theorem~\ref{intro-th.1}
yields also
\be{old.6}
\lim_{\Lambda\nearrow\Z^d}\lim_{n\to\infty} \frac{1}{\sqrt n} \log P\pare{
||\ind_{\Lambda} l_\infty||_2^2> n}=-\I(2).
\ee
Indeed, in order to use our subadditive argument, Lemma~\ref{lem-sub.1},
we need first to show that for some $\gamma>0$, for any
$\alpha$ large enough, and for $n$ large enough
\be{lawler.10bis}
\begin{split}
P_0\big(||\ind_\Lambda&l_{\lfloor\alpha\sqrt{n}\rfloor}
||_2^2\ge n,\ S_{\lfloor\alpha\sqrt{n}\rfloor}=0\big)\\
&\le P\pare{ ||\ind_{\Lambda} l_\infty||_2^2> n}
\le n^\gamma P_0\pare{||\ind_\Lambda l_{\lfloor\alpha\sqrt{n}\rfloor}
||_2^2\ge n,\ S_{\lfloor\alpha\sqrt{n}\rfloor}=0}.
\end{split}
\ee
The upper bound in \reff{lawler.10bis}
is obtained from Proposition~\ref{prop-time.1},
whereas the lower bound is immediate.

Now, we proceed with the link with intersection local times.
First, as mentioned in \reff{CM-finite},
Chen and M\"orters prove also that for any finite $\Lambda\subset \Z^d$
\[
\lim_{n\to\infty}\frac{1}{n^{1/2}}\log P\pare{ \bra{\ind_{\Lambda} l_\infty,
\tilde l_\infty}>n}=-2I_{CM}(\Lambda),
\]
with $I_{CM}(\Lambda)$ converging to $I_{CM}$ as $\Lambda$ increases to cover $\Z^d$.
The important feature is that for any fixed $\epsilon>0$,
we can fix a finite $\Lambda$ subset of $\Z^d$ such that $|I_{CM}(\Lambda)-I_{CM}|\le 
\epsilon$.
Note now that by Cauchy-Schwarz' inequality, and for finite set $\Lambda$
\be{old.7}
\bra{\ind_{\Lambda}l_\infty,\tilde l_\infty}\le 
||\ind_{\Lambda}l_\infty||_2 ||\ind_{\Lambda}\tilde l_\infty||_2.
\ee
Inequalities \reff{old.6} and \reff{old.7} imply by routine consideration that
\be{old.8}
\limsup_{\Lambda\nearrow\Z^d} \limsup_{n\to\infty}
\frac{1}{\sqrt n} \log P\pare{\bra{1_\Lambda l_{\infty},\tilde l_\infty}>n}
\le -\I(2)\inf_{\alpha>0}\acc{{\sqrt \alpha}+\frac{1}{{\sqrt \alpha}}}=-2\I(2).
\ee
When ${\bf k_n^*}$ is the sequence which enters into defining
$\A_n^*(1,\Lambda)$ in \reff{lawler.3} (see also \reff{lawler.2}),
we have the lower bound 
\be{old.9}
P\pare{ \bra{\ind_{\Lambda}l_\infty,\tilde l_\infty}>n}
\ge  P\pare{ l_{\lfloor \alpha\sqrt n\rfloor}|_{\Lambda}=
{\bf k^*_n}|_{\Lambda}, S_{\lfloor \alpha\sqrt n\rfloor}=0}^2.
\ee
Following the same argument as in the proof of Section~\ref{sec-ergodic}, 
we have
\be{old.10}
\liminf_{n\to\infty}  \frac{1}{\sqrt n} \log P\pare{\bra{l_{\infty},\tilde l_\infty}>n}
\ge -2\I(2).
\ee
\reff{old.8} and \reff{old.10} conclude the proof \reff{identification}.

\section{Applications to RWRS.}
\label{sec-rwrs}
We consider a certain range of parameters $\acc{(\alpha,\beta): 
1<\alpha<\frac{d}{2},
1-\frac{1}{\alpha+2}<\beta<1+\frac{1}{\alpha}}$, which we have called Region II in
\cite{AC05}.  Also, if $\Gamma(x)=\log(E[\exp(x \eta(0))])$, then there
are positive constants $\Gamma_0$ and $\Gamma_{\infty}$ (see \cite{AC05}), such that
\be{rwrs.2}
\lim_{x\to 0} \frac{\Gamma(x)}{x^2}=\Gamma_0,\quad\text{and}\quad
\lim_{x\to \infty} \frac{\Gamma(x)}{x^{\alpha^*}}=\Gamma_{\infty},\quad\text{and}\quad
\frac{1}{\alpha}+\frac{1}{\alpha^*}=1.
\ee
A classical way of obtaining large deviations is through exponential bounds for
$\P(\bra{\eta,l_n}\ge y n^\beta)$. For instance, if we expect the latter quantity
to be of order $\exp(-c n^\zeta)$, then a first tentative would be to optimize
over $\lambda>0$ with $b=\beta-\zeta$ in the following
\ba{tent.1}
\P\pare{\bra{\eta,l_n}\ge y n^\beta}&\le & e^{-\lambda n^{\beta-b}}E[
\exp\pare{\lambda \frac{\bra{\eta,l_n}}{n^b}}]\cr
&\le & 
e^{-\lambda n^{\zeta}}
E_0[\exp\pare{\sum_{z\in \Z^d} \Gamma\pare{\frac{\lambda l_n(z)}{n^b}}}].
\ea
We need to distinguish
asymptotic regimes at zero or at infinity for $\Gamma(\frac{\lambda l_n(z)}{n^b})$ 
according to whether $l_n(x)<n^{b-\epsilon}$ or $l_n(x)> n^{b+\epsilon}$ respectively.
For $\epsilon>0$, we introduce
\[
\Dh=\acc{x \in \Z^d:\ l_n(x) \geq n^{b+\epsilon}},\qquad
\Db= \acc{x \in \Z^d:\  0<l_n(x) \leq n^{b-\epsilon}},
\]
and,
\[
\RR_{\epsilon}=\acc{x \in \Z^d; n^{b-\epsilon}
\le l_n(x) \le n^{b+\epsilon}}.
\]
Then, for any $\epsilon_0>0$ small
\be{rwrs.3}
\P\pare{\bra{\eta,l_n}\ge yn^{\beta}}\le 
\P\pare{\bra{\eta,\ind_{\Dh}l_n}\ge 
(1-\epsilon_0)yn^{\beta}}+I_1+I_2,
\ee
where
\be{def-I}
I_1:=\P\pare{\bra{\eta,\ind_{\Db}l_n}\ge \frac{\epsilon_0}{2} yn^{\beta}}
,\quad\text{and}\quad
I_2:=\P\pare{\bra{\eta,\ind_{\RR_{\epsilon}}l_n}\ge \frac{\epsilon_0}{2} yn^{\beta}} .
\ee
We have now to show that the contribution of $\Db$ and $\RR_{\epsilon}$ which
concerns the {\it low} level sets, is negligible. We gather the two estimates
in the next subsection. We treat afterwards $\Dh$.

\subsection{Contribution of {\it small} local times.}
We first show that $I_1$ is negligible. 
Set $\B=\acc{||\ind_{\Db}l_n||_2^2\ge \delta n^{\beta+b}}$, 
for a $\delta>0$ to be chosen later. For any $\lambda>0$
\be{rwrs.4}
\P\pare{\bra{\eta,\ind_{\Db}l_n}\ge \frac{\epsilon_0}{2} yn^{\beta}}\le 
P(\B)+
e^{-\lambda n^{\beta-b} \frac{\epsilon_0}{2} y}
E_0\cro{\ind_{\B^c} \exp\pare{ \sum_{\Db} \Gamma\pare{\frac{\lambda l_n(x)}{n^b}}}}.
\ee
Now, for any $\lambda>0$ and $n$ large enough, we have for $x\in \Db$ that
\[
\Gamma(\frac{\lambda l_n(x)}{n^b})\le \Gamma_0(1+\epsilon_0) (\frac{\lambda l_n(x)}{n^b})^2,
\]
so that
\be{rwrs.5}
\P\pare{\bra{\eta,\ind_{\Db}l_n}\ge \frac{\epsilon_0}{2} yn^{\beta}}\le P(\B)+
\exp\pare{-n^\zeta\pare{ \lambda 
\frac{\epsilon_0}{2} y- \lambda^2 \Gamma_0(1+\epsilon_0)\delta}}.
\ee
Since $\beta+b>1$, Lemma 1.8 of \cite{AC05} gives that $-\log\pare{P(\B)}\ge M n^{\zeta}$,
for any $\delta>0$, and any large constant $M$.
Finally, for any $\epsilon_0$ fixed, and a large constant $M$,
we first choose $\lambda$ so that $\lambda 
\frac{\epsilon_0}{2} y\ge 2M$. Then, we choose $\delta$ small enough
so that $\lambda \Gamma_0\delta\le \frac{\epsilon_0}{4} y$.

We consider the contribution of $\RR_{\epsilon}$. We use here our hypothesis
that the $\eta$ are bell-shaped random variables, since it leads to clearer
derivations. Thus, according to Lemma 2.1 of~\cite{AC04}, we have
\be{rwrs.6}
\P\pare{ \bra{\eta,1\{\RR_{\epsilon}\}l_n}\ge yn^\beta}\le 
\P\pare{\sum_{\RR_{\epsilon}}\eta(x)\ge n^{\beta-b-\epsilon}}.
\ee
By Proposition 1.9 of \cite{AC05}, we can assume that $|\RR_{\epsilon}|< n^\gamma$, with
\be{rwrs.7}
\gamma<\gamma_0:=\frac{1}{1-\frac{2}{d}}\frac{\alpha-1}{\alpha+1}\beta=
\frac{1-\frac{1}{\alpha}} {1-\frac{2}{d}}\zeta<\zeta\quad\text{if}\quad \alpha<\frac{d}{2}.
\ee
Note that $\gamma_0$ given in \reff{rwrs.7} is lower than $\zeta$ when $\alpha<d/2$. Using Lemma A.4 of \cite{AC04}, we obtain
\be{rwrs.8}
\P\pare{\sum_{\RR_{\epsilon}} \eta(x)\ge n^{\beta-b-\epsilon}, |\RR_{\epsilon}|\le
n^\gamma}\le \exp\pare{ -Cn^{\gamma+\alpha(\beta-b-\epsilon-\gamma)}}.
\ee
For the left hand side of \reff{rwrs.8} to be negligible, we would need
(recall that $\alpha>1$)
\be{rwrs.9}
\gamma+\alpha(\beta-b-\epsilon-\gamma)>\beta-b\Longleftrightarrow
(\beta-b)(\alpha-1)>(\alpha-1)\gamma\Longleftrightarrow
\zeta>\gamma.
\ee
This last inequality has already been noticed to hold in \reff{rwrs.7}.

\subsection{Contribution of {\it large} local times.}
\subsubsection{Upper Bound}
We deal now with the contributions of $\Dh$. For any $\lambda>0$ (recalling that $\beta
-b=\zeta=\frac{\alpha\beta}{\alpha+1}$)
\be{rwrs.10}
\P\pare{\bra{\eta,1\{\Dh\}l_n}\ge (1-\epsilon_0) y n^\beta}\le 
e^{-\lambda n^{\beta -b}(1-\epsilon_0)y}
E_0\cro{ \exp \pare{\sum_{x \in \Dh} \Gamma\pare{ \frac{\lambda l_n(x)}{n^{b}}}}}.
\ee
Now, for $\lambda$ {\it not too small}, when $n$ is large enough we have 
\be{rwrs.11}
\sum_{x \in \Dh} \Gamma\pare{\frac{\lambda l_n(x)}{n^{b}}}\le(\Gamma_{\infty}+\epsilon_0)
\lambda^{\alpha^*}\pare{\frac{||\ind_{\Dh}l_n||_{\alpha^*}}{n^{b}}}^{\alpha^*}.
\ee
Thus, \reff{rwrs.10} becomes
\be{rwrs.12}
\P\pare{\bra{\eta,1\{\Dh\}l_n}\ge (1-\epsilon_0) y n^\beta}\le
\exp\pare{ -n^\zeta\pare{ \lambda (1-\epsilon_0)y-\lambda^{\alpha^*}(\Gamma_{\infty}+\epsilon_0)
\frac{||\ind_{\Dh}l_n||_{\alpha^*}^{\alpha^*}}{n^{b\alpha^*+\zeta}}}}.
\ee
Now, optimizing in $\lambda$ in the right hand side of \reff{rwrs.12}, we obtain
\be{rwrs.13}
(1-\epsilon_0)y=\alpha^* (\Gamma_{\infty}+\epsilon_0)
\frac{||\ind_{\Dh}l_n||_{\alpha^*}^{\alpha^*}}{n^{b\alpha^*+\zeta}}\lambda^{\alpha^* -1}.
\ee
Now, recall that in order to fall in the asymptotic regime of $\Gamma$ at infinity, we assumed
that $\lambda$ were not too small. In other words, 
in view of \reff{rwrs.13}, we would
need a bound of the type $||\ind_{\Dh}l_n||_{\alpha^*}\le An^{\zeta}$ for a large constant $A$.
Now, using 
Proposition~\ref{prop-alpha}, there is a constant $\I(\alpha^*)$ such that 
\be{rwrs.14}
\lim_{n\to\infty} \frac{1}{n^\zeta} \log
P\pare{||\ind_{\Dh}l_n||_{\alpha^*}\ge A n^{\zeta}}\le -\I(\alpha^*) A .
\ee
Thus, we can assume that $\lambda$ satisfying \reff{rwrs.13} is bounded from below.
Also, replacing the value of $\lambda$ obtained in \reff{rwrs.13} in inequality \reff{rwrs.12},
and using that $\Gamma_{\infty}^{-1}=\alpha^*(\alpha c_{\alpha})^{\alpha^*-1}$, we find that
\be{rwrs.15}
\P\pare{\bra{\eta,\ind_{\{\Dh\}}l_n}\ge (1-\epsilon_0) y n^\beta}\le
E_0\cro{\exp\pare{ -c_{\alpha} (1-\delta_0)
\pare{ \frac{yn^\beta}{||\ind_{\Dh}l_n||_{\alpha^*}}}^\alpha}},
\ee
where $(1-\delta_0)=(1-\epsilon_0)^\alpha(1+\epsilon_0)^{1-\alpha}$, 
which can be made as close
as 1, as one wishes. Now, it is easy to conclude that
\ba{upper-LDP}
\limsup_{n\to\infty} \frac{1}{n^\zeta} \log 
P\pare{\bra{\eta,l_n}\ge  y n^\beta}&\le&
-c_{\alpha} \inf_{\xi>0}\acc{ \pare{\frac{y}{\xi}}^{\alpha} + \I(\alpha^*) \xi}\cr
&=& -c_{\alpha} (\alpha+1) \pare{ \frac{y \I(\alpha^*) }{\alpha}}^{\frac{\alpha}{\alpha+1}}.
\ea

\subsubsection{Lower Bound for RWRS.}
We call in this section $\bar\D=\acc{z\in \Z^d: l_n(z)\ge \delta n^\zeta}$, for a fixed
but small $\delta$. Since, we have assumed the $\eta$-variables to have a bell-shaped
distribution, we have according to Lemma 2.1 of~\cite{AC04},
\be{rwrs.16}
\P\pare{\bra{\eta,l_n}\ge y n^\beta}\ge \P\pare{\bra{\eta,\ind_{\bar\D}l_n}\ge y n^\beta}.
\ee
Then, we condition on the random walk law, and average with respect to the $\eta$ variables
which we require to be large on each site of $\bar \D$. Recall now that we can
assume $|\bar\D|\le 1/\delta^2$ by \reff{level.21} (for $\delta$ small enough). We use \reff{rwrs.16} to deduce for any $\epsilon>0$
\[
\begin{split}
\P(\bra{\eta,l_n}\ge y n^\beta)\ge&E_0
\cro{\Q\cro{\min_{z\in \bar\D} \eta(z)\ge \epsilon n^\zeta,\ 
\bra{\eta,\ind_{\bar\D}l_n}\ge y n^\beta}}\\
\ge & E_0\cro{\sup_{x(i),i\in \bar\D}
\acc{C^{|\bar\D|}\exp\pare{-c_{\alpha}(1+\epsilon)
\sum_{i\in \bar\D} x(i)^{\alpha}}: \bra{x,\ind_{\bar\D}l_n}\ge yn^\beta}}\\
\ge & C^{1/\delta^2} E_0\cro{\ind\acc{|\D|\le \frac{1}{\delta^2}}
\exp\pare{-c_{\alpha}(1+\epsilon) 
\pare{\frac{yn^\beta}{||\ind_{\bar\D}l_n||_{\alpha^*}}}^{\alpha}}}\\
\ge & C^{1/\delta^2}  \exp\pare{-c_{\alpha}(1+\epsilon)
\pare{\frac{yn^\beta}{\xi^* n^\zeta}}^{\alpha}}
P\pare{ ||\ind_{\bar\D}l_n||_{\alpha^*}\ge \xi^* n^\zeta,\ |\D|\le \frac{1}{\delta^2}},
\end{split}
\]
where $\xi^*$ realizes the infimum in \reff{upper-LDP}. Now,
as $\epsilon$ is sent to 0 after $n$ is sent to infinity, we obtain
\be{lower-LDP}
\liminf_{n\to\infty} \frac{1}{n^\zeta} \log
P\pare{\bra{\eta,l_n}\ge  y n^\beta}\ge
-c_{\alpha} (\alpha+1) \pare{ \frac{y \I(\alpha^*) }
{\alpha}}^{\frac{\alpha}{\alpha+1}}.
\ee

\section{Appendix.}\label{sec-appendix}
%PROOF OF THE LEMMA.
\subsection{Proof of Lemma~\ref{circuit-lem.1}.}
Fix ${\bf k}\in V(\Lambda',n)$. By Chebychev's inequality, for any $\lambda>0$
\be{circuit.15}
P\big( \sum_{i=1}^{|{\bf k}|}
1_{\acc{|S_{T^{(i)}}-S_{T^{(i-1)}}|>{\sqrt L}, T^{(|{\bf k}|)}<\infty}}\ge
\epsilon\sqrt{n}\big) \le e^{-\lambda \epsilon \sqrt{n}}
E\cro{\prod_{i=1}^{|{\bf k}|}e^{\lambda \ind\acc{|S_{T^{(i)}}-S_{T^{(i-1)}}|>{\sqrt L}}}}.
\ee
Now, by using the strong Markov's property, 
and induction, we bound the right hand side of
\reff{circuit.15} by
\be{circuit.16}
e^{-\lambda \epsilon \sqrt{n}} \pare{ \sup_{z\in \Lambda'\cup\{0\}} 
E_{z}\cro{ e^{\lambda \ind\acc{|S_T|>\sqrt{L}, T<\infty}}}}^{|{\bf k}|}.
\ee
Now, 
\ba{circuit.17}
E_{z}\cro{ e^{\lambda \ind\acc{|S_T|>\sqrt{L},T<\infty}}}&\le&
1+(e^\lambda-1) P_{z}(|S_T|>\sqrt{L},T<\infty)\cr
&\le & 1+(e^\lambda-1) P_{z}\pare{ \cup\acc{T(\xi)<\infty:\ \xi\in \Lambda;\ |\xi-z|>\sqrt{L}}}\cr
&\le & 1+(e^\lambda-1)|\Lambda| \sup\acc{P_z(T(\xi)<\infty);\ |\xi-z|>\sqrt{L}, \xi\in \Lambda}\cr
&\le & 1+(e^\lambda-1) \frac{\bar c|\Lambda|}{L^{d/2-1}}
\le \exp\pare{ (e^\lambda-1)\frac{\bar c|\Lambda|}{L^{d/2-1}}}.
\ea
Now, since $|{\bf k}|\le c_0 \sqrt{n}$, we have
\be{circuit.18}
P\pare{ \sum_{i=1}^{|{\bf k}|} 
1_{\acc{|S_{T^{(i)}}-S_{T^{(i-1)}}|>{\sqrt L}, T^{(|{\bf k}|)}<\infty}}\ge
\epsilon\sqrt{n}}\le \exp\pare{-\sqrt{n}
\pare{\lambda\epsilon- (e^\lambda-1)\frac{c_0\bar c|\Lambda|}{L^{d/2-1}}}}.
\ee
Thus, for any $\epsilon>0$, we can choose $L$ large enough so that the result holds.
\qed
\subsection{Proof of Lemma~\ref{lem-encageS}.}
We first introduce a fixed scale, $l_0\in \N$, 
to be adjusted later as a function of
$|\Lambda|$, and assume that $|x-y|\ge 4 |\Lambda| l_0$. 
Indeed, the case $|x-y|\le 4 |\Lambda|l_0$ is easy to treat since $P_x(S_T=y)>0$ implies
the existence of a path from $x$ to $y$ avoiding $\Lambda$; it is then easy to see
that since $\Lambda$ is finite, the length of the shortest path
joining $x$ and $y$ and avoiding $\Lambda$ can be bounded by a constant
depending only on $|\Lambda|$. Forcing the walk to follow this path costs
only a positive constant which depends on $|\Lambda|$.

We introduce two sets of concentric shells around $x$ and $y$: for $i=1,\dots,|\Lambda|-1$
\be{encage.3}
C_i=B(x,(2i+2)l_0)\bs B(x,2il_0),\quad\text{and}\quad C_0=B(x,2l_0),
\ee
and similarly $\acc{ D_i,i=0,\dots,|\Lambda|}$ 
are centered around $y$, and for all
$i,j$ $C_i\cap D_j=\emptyset$. There is necessarely $i,j\le |\Lambda|$ 
such that
\be{encage.4}
C_i\cap \Lambda=\emptyset,\quad\text{and}\quad D_j\cap \Lambda=\emptyset.
\ee
Define now two stopping times corresponding to exiting {\it mid-}$C_i$ and
entering {\it mid}-$D_j$
\be{encage.5}
\sigma_i=\inf\acc{n\ge 0:\ S_n\not\in B(x,(2i+1)l_0)},\quad\text{and}\quad 
\tau_j=\inf\acc{n\ge 0:\ S_n\in B(y,(2j+1)l_0)}.
\ee
Note that when $\sigma_i<\infty$ and $\tau_j<\infty$, we have
$\dist(S_{\sigma_i},\Lambda)\ge l_0$, and $\dist(S_{\tau_j},\Lambda)\ge l_0$.
We show that for any $L$ we can find $\epsilon_L$ 
(going to 0 as $L\to\infty$), such that
\be{encageS-key}
P_x(T(\S)<T<\infty,\ S_T=y)\le 
\frac{\e_L}{2} P_x(S_T=y).
\ee
Note that \reff{encageS-key} implies that for $\e_L$ small enough
\be{encage.15}
P_x(S_T=y)\le \frac{1}{1-\e_L/2} P_x(T<T(\S),\ S_T=y)\le
e^{\e_L} P_x(T<T(\S),\ S_T=y).
\ee
To show \reff{encageS-key},
we condition the flight $\acc{S_0=x,S_T=y}$ 
on its values at $\sigma_i$ and $\tau_j$
\be{encage.6}
P_x(S_T=y)\ge \sum_{z\in C_i} P_x\pare{S_{\sigma_i}=z,\ \sigma_i<T}P_z(\tau_j<T)
\inf_{z'\in D_j}P_{z'}(S_T=y).
\ee
Note that if $P_x(S_T=y)>0$, there is necessarely a path from $D_j$ to $y$ which avoids
$\Lambda$ so that, there is a constant $c_0$ (depending only on $l_0$) such that
\be{encage.7}
\inf_{z'\in D_j}P_{z'}(S_T=y)>c_0.
\ee
We need to estimate $P_z(\tau_j<T)$. First, 
by classical estimates (see Proposition 2.2.2 of \cite{LAWLER}), 
there are $c_1,c_2>0$ such that
when $|x-y|\ge 4l_0|\Lambda|$, and $z\in C_i$
\be{encage.8}
\frac{c_2\capa(D_j)}{|z-y|^{d-2}}\le P_{z}(\tau_j<\infty)\le 
\frac{c_1\capa(D_j)}{|z-y|^{d-2}}.
\ee
We establish now that if we choose $l_0$ so that
\be{encage.11}
l_0^{d-2}\ge 2 |\Lambda| \frac{c_1c_G2^{d-2}}{c_2},
\quad\text{then}\quad P_z(\tau_j<T)\ge \frac{1}{2} P_z(\tau_j<\infty). 
\ee
Since $\dist(z,\Lambda)>l_0$
\ba{encage.9}
P_{z}(T<\tau_j<\infty)&\le &\sum_{\xi\in \Lambda\bs D_0} 
P_{z}\pare{S_T=\xi,\ T<\tau_j<\infty} P_{\xi}(\tau_j<\infty)\cr
&\le & |\Lambda| \sup_{\xi\in \Lambda\bs D_0}
\acc{P_z(T(\xi)<\infty)P_{\xi}(\tau_j<\infty)}.
\ea
We use again estimate \reff{encage.8} to obtain
\be{encage.91}
P_{z}(T<\tau_j<\infty)\le c_1 c_G |\Lambda| \sup_{\xi \in \Lambda\bs D_0} \acc{
\frac{1}{|z-\xi|^{d-2}}\times \frac{\capa(D_j)}{|\xi -y|^{d-2}}}.
\ee
Now, for $\xi \in \Lambda\bs D_0$, we have $\min(|z-\xi|,|\xi-y|)>l_0$, 
and on the other side
the triangle inequality yields $\max(|z-\xi|,|\xi-y|)> \frac{|z-y|}{2}$.
Thus, we obtain
\ba{encage.10}
P_{z}(T<\tau_j)&\le &\frac{c_1c_G2^{d-2}}{l_0^{d-2}} |\Lambda| \frac{\capa(D_j)}{|z-y|^{d-2}}\cr
&\le & \frac{c_1c_G2^{d-2}}{c_2}\frac{|\Lambda|}{l_0^{d-2}}\frac{c_2\capa(D_j)}{|z-y|^{d-2}}\cr
&\le & \frac{c_1c_G2^{d-2}}{c_2}\frac{|\Lambda|}{l_0^{d-2}} P_z(\tau_j<\infty).
\ea
This implies \reff{encage.11}.

Now, for any $z\in C_i$, by conditioning
on $S_{T(\S)}$, we obtain
\be{encage.12}
P_z(T(\S)<T<\infty,\ S_T=y)\le E_z\cro{ \ind\acc{T(\S)<T<\infty}P_{S_{T(\S)}}
\pare{S_T=y}}\le \frac{c_G}{L^{d-2}}.
\ee
Thus, for any $z\in C_i$,
\be{encage.13}
P_z(\tau_j<T) \inf_{z'\in D_j}P_{z'}(S_T=y)\ge c_0\frac{c_2\capa(D_j)}{2|z-y|^{d-2}}\ge
\frac{P_z(T(\S)<T<\infty,\ S_T=y)}{\e_L/2},
\ee
with (recalling that $|x-y|\ge 4|\Lambda|l_0$ and $|x-y|\le \sqrt{L}$), with a constant
$C(\Lambda)>0$
\be{encage.14}
\e_L=\frac{4 c_G |z-y|^{d-2}}{c_0c_2 \capa(D_j) L^{d-2}}\le 
\frac{ 2^d c_G}{c_0c_2 \capa(D_j)}\pare{ \frac{|x-y|}{L}}^{d-2}\le
C(\Lambda) \pare{ \frac{1}{\sqrt{L}}}^{d-2}.
\ee
Now, after summing over $z\in C_i$, we obtain \reff{encageS-key}.
\qed
\subsection{Proof of Lemma~\ref{lem-encageB}.}
We consider two cases: (i) $\sqrt{L}< |x-y|\le \kappa L$ where $\kappa$ is 
a small parameter, and (ii) $|x-y|> \kappa L$. 

Also, we denote by $C(\lambda)$ a positive constant which depend only
on $|\Lambda|$. We might use the same name in different places.

\noindent{\underline{Case (i).}} We use the same steps
as in the previous proof up to \reff{encage.13} where we replace $|z-y|$ by $2|x-y|$,
and obtain
\be{encage.16}
P_x(S_T=y)\ge \frac{c_0}{2} \frac{c_2\capa(D_j)}{2^{d-2}|x-y|^{d-2}}.
\ee
Now, \reff{encage.12} implies that if
\be{encage.17}
\kappa^{d-2}\le \frac{c_0c_2 \capa(D_j)}{2^{d}c_G},\quad\text{then}\quad
P_x(S_T=y)\le 2 P_x(T<T(\S),\ S_T=y).
\ee
\noindent{\underline{Case (ii).}} First note that
\be{encage.fact1}
P_x(S_T=y)\le P_x(T(y)<\infty) \le \frac{c_G}{|x-y|^{d-2}}.
\ee
Now, set $L'=\kappa L$, and note that
$\diam(\C)$ is a multiple (depending only on $\Lambda$) times $L'$. Now,
a way of realizing $\acc{ S_T=y, T<T(\S)}$ is to go through a finite number 
of adjacent spheres of diameter $L'$. From a hitting point on one sphere, 
we force the walk to exit
only from a tiny fraction of the surface of the next sphere, until we reach the last sphere, say
on $z^*$, for which it is easy to show that there are two universal positive constants
$c,c'$ such that
\be{encage.fact2}
P_{z^*}(S_T=y, T<T(\S))\ge c P_{z^*}(T(y)<\infty)\ge c' \frac{\tilde c_G}{|x-y|^{d-2}}.
\ee

Note that when starting on $x$, the
probability of exiting $B(x,|x-y|)$ through site $y$ is of order of the surface $|x-y|^{1-d}$,
and this is much smaller of $P_x(T(y)<\infty)$ which should be close to $P_x(S_T=y)$ in
cases where all other points of $\Lambda$ be very far from $x,y$. 
Thus, we have to consider
more paths than $\acc{S_{T(B(x,|x-y|)^c)}=y,S_0=x}$. 
By Lemma~\ref{lem-cluster} and Remark~\ref{rem-cluster},
there is a finite sequence $x_1,\dots,x_k$
(not necessarely in $\C$) 
such that $L'/2\le |x_{i+1}-x_i|\le L'$ and such that
$B(x_i,L)\subset \S(\C)$.
\be{encage.18}
\delta=\frac{1}{4|\Lambda|^{\frac{1}{d-1}}},\quad
Q_i=\acc{z: |z-x_i|=|x_{i+1}-x_i|},\quad\text{and}\quad
\Sigma_i=Q_i\cap B(x_{i+1},\frac{L'}{4}).
\ee
Note that $|\Sigma_i|$ is of order $(\frac{L'}{4})^{d-1}$. We can throw 
$|\Lambda|$ points
on $\Sigma_i$, say at a distance of at least $\delta L'$,
and one of them, say $y_{i}^*$, necessarely satisfies 
\be{encage.19}
B(y_{i}^*,\delta L')\cap \Lambda =\emptyset,\quad\text{and set}\quad
B_i^*=B(y_{i}^*,\frac{\delta L'}{2})\cap \Sigma_i.
\ee
Now, when the walk starts on $x_{i+1}$, it exits from any point $z\in Q_{i+1}$ with
roughly the same chances (see i.e. Lemma 1.7.4 of~\cite{LAWLER}), 
so that there is $c_S$ such that for $i\ge 0$,
\be{encage.20}
P_{x_{i+1}}(S_{H_i}=z)\ge \frac{c_S}{|x_{i+2}-x_{i+1}|^{d-1}} ,\quad\text{where}\quad
H_i:=T(Q_{i+1}).
\ee
By Harnack's inequality (see Theorem 1.7.2 of~\cite{LAWLER}), 
for any $z\in B_i^*$
\be{encage.21}
P_z\pare{ S_{H_i}\in B^*_{i+1}}\ge \frac{c_S |B^*_{i+1}|}{(2L)^{d-1}}
\ee
Now, there is $\chi>0$ such that
\[
|B^*_i|\ge \chi(\delta L')^{d-1},
\]
which yields
\be{encage.key}
P_z\pare{ S_{H_i}\in B^*_{i+1}}\ge c_S \chi (\frac{\delta \kappa}{2})^{d-1}.
\ee
Note that it costs more to hit $\Lambda$ before $Q_{i+1}^c$. Indeed,
\ba{encage.22}
P_z\pare{S_{H_i}\in B^*_{i+1},\ T<H_i}&\le& 
\sum_{\xi\in \Lambda} P_z(T(\xi)<\infty) P_{\xi}\pare{ H_i<\infty}\cr
&\le &\sup_{\xi \in \Lambda} \frac{c_G|\Lambda|}{|z-\xi|^{d-2}}\times
\frac{c_1 \capa\pare{ B^*_{i+1}}}{|\xi-y_{i+1}^*|^{d-2}}.
\ea
By definition, $
\capa( B^*_{i+1} )\le |B^*_{i+1}|\le \chi (\delta L')^{d-1}$.
Now, $z$ and $y_{i+1}^*$ are chosen in such a way that 
$\min(|z-\xi|,|\xi-y_{i+1}^*|)\ge \frac{\delta L'}{2}$ so that
\be{encage.23}
P_z\pare{S_{H_i}\in B^*_{i+1},\ T<H_i}\le
\frac{ c_Sc_1 \chi |\Lambda| (\delta L')^{d-1} }
{ (\delta L'/2)^{2d-4}}.
\ee
Since in $d\ge 5$, we have $2d-4>d-1$, $L$ can be chosen large enough so that 
\be{encage.24}
P_z\pare{S_{H_i}\in B^*_{i+1}, \ H_i<T}\ge \frac{1}{2} P_z\pare{ S_{H_i}\in B^*_{i+1}}
\ee
Now, we define $\theta_k$ as the time-translation of $k$ units of a random walk trajectory,
and $\tilde H_i=H_i\circ \theta_{H_{i-1}}$. The following scenario produces
$\acc{ S_T=y,\ T<T(\S)}$: 
\be{encage.25}
\bigcap_{i=1}^k \acc{S_{\tilde H_i}\in B^*_{i+1},\  \tilde H_i <T\circ \theta_{H_{i-1}}} 
\cap \acc{S_{T\circ \theta_{H_k}}=y,T\circ\theta_{H_k}<T(\S)
\circ\theta_{H_k}}
\ee
By using the strong Markov's property, and \reff{encage.24}, we obtain
\be{encage.26}
P_x\pare{S_T=y,\ T<T(\S)}\ge \pare{ \frac{c_S\chi}{2} (\frac{\delta \kappa}{2})^{d-1}}^k
\inf_{z\in B^*_{k}} P_z\pare{S_T=y,T<T(\S)}.
\ee
In the last term in \reff{encage.26}, note that for any $z\in B_k^*$,
$L'/2\le |z-y|\le L'$ so that we are in the situation of Case(i),
where inequality \reff{encage.17}, and \reff{encage.16} yields
\[
P_z\pare{S_T=y,T<T(\S)}\ge \frac{1}{2} P_z\pare{S_T=y}\ge \frac{c}{|z-y|^{d-2}}.
\]
Since Lemma~\ref{lem-cluster} establishes that for some
constant $C(\Lambda)>0$, $\diam(\C)\le C(\Lambda)L$, and
$|z-y|\ge \frac{\kappa}{2} L$, we have for a constant $C(\Lambda)$ 
\[
P_x\pare{S_T=y,\ T<T(\S)}\ge \frac{C(\Lambda)}{|x-y|^{d-2}}\ge
\frac{C(\Lambda)}{c_G} P_x\pare{S_T=y}.
\]
\qed
\subsection{Proof of Lemma~\ref{trans-lem.1}.}
We start with shorthand notations $\S_1=\S(\C)$ and $\tilde
\S_1=\S(\T(\C))$, and we define
\[
\S_2=\acc{z: \dist(z,\C)=2\max(\diam(\C),L)},
\]
and $\tilde \S_2$ is similar to $\S_2$ but
$\T(\C)$ is used instead of $\C$ in its definition.

First, we obtain an 
upper bound for the weights of paths joining $y$ to $x$ by
conditioning over hitting sites on $\S_2$ and $\S_1$, and by using
the strong Markov's property
\ba{trans.3}
P_y(S_T=x)&=& 
\sum_{z_1\in \S_1} E_y\cro{ \ind_{\acc{T(\S_2)<T}}
P_{S_{T(\S_2)}}(S_{T(\S_1)}=z_1,T(\S_1)<T)}P_{z_1}(S_T=x)\cr
&\le & P_y(T(\S_2)<\infty)\sum_{z_1\in \S_1} 
\pare{ \sup_{z\in \S_2}P_z(S_{T(\S_1)}=z_1)} P_{z_1}(S_T=x)
\ea
We need to compare \reff{trans.3} 
with the corresponding decomposition for
trajectories starting on $y$ with $\acc{S_T=\tilde x}$, where
we set $\tilde x=\T(x)$ for simplicity,
\be{trans.15}
\begin{split}
P_y(S_T=\tilde x)=& 
\sum_{\tilde z_1\in \tilde \S_1} E_y\cro{ \ind_{\acc{T(\tilde \S_2)<T}}
P_{S_{T(\tilde \S_2)}}(S_{T(\tilde \S_1)}=\tilde z_1,T(\tilde \S_1)<T)}
P_{\tilde z_1}(S_T=\tilde x)\\
\ge & P_y(T(\tilde\S_2)<T)\!\!\sum_{\tilde z_1\in \tilde \S_1} 
\!\!\inf_{\tilde z\in \tilde \S_2} P_{\tilde z} 
\pare{S_{T(\tilde \S_1)}=\tilde z_1,
T(\tilde \S_1)<T} P_{\tilde z_1}(S_T=\tilde x).
\end{split}
\ee
We now bound each term in \reff{trans.3} by 
the corresponding one in \reff{trans.15}.

\underline{ About $P_{z_1}(S_T=x)$}.
From \reff{clus.3} of Lemma~\ref{lem-cluster}, 
$\S_2\cap\Lambda=\C$.
By the same reasoning as in the proof of
Lemma~\ref{lem-encageB}, 
there is a constant $C_0$ such that for any $z_1\in \S_1$
\be{trans.4}
P_{z_1}(S_T=x)\le C_0 P_{z_1}(S_T=x,T<T(\S_2)).
\ee
As long as we consider paths from $\S_1$ to $x$ 
which do not escape $\S_2$, we can
transport them, using translation invariance of the law of random walk 
\be{trans.5}
P_{\tilde z_1}(S_T=\tilde x,T<T(\tilde \S_2))=P_{z_1}(S_T=x,T<T(\S_2)),
\ee
and by using \reff{trans.4} and \reff{trans.5}, we finally obtain
\be{trans.6}
P_{z_1}(S_T=x)\le C_0 P_{\tilde z_1}(S_T=\tilde x,T<T(\tilde \S_2))
\le C_0 P_{\tilde z_1}(S_T=\tilde x).
\ee
\underline{ About $P_y(T(\S_2)<\infty)$}. 
By Proposition 2.2.2 of~\cite{LAWLER}, there are $c_1,c_2$ positive constants such that
\be{trans.7}
\frac{c_2 \capa(\S_2)}{|y-x|^{d-2}}\le P_y(T(\S_2)<\infty)\le 
\frac{c_1 \capa(\S_2)}{|y-x|^{d-2}},
\ee
and \reff{trans.7} holds also
with a tilda over $x$ and $\S_2$. Since $|y-\tilde x|\le
2 |y-x|$ by \reff{clus.7}, we have
\be{trans.8}
P_{y}(T(\S_2)<\infty)\le \frac{c_1}{c_2} 2^{d-2} P_{y}(T(\tilde \S_2)<\infty).
\ee
We need now to check that paths reaching $\tilde \S_2$ from $y$ have 
good chances not to meet any sites of $\Lambda$. In other words, we need
\be{trans.40}
P_y(T(\tilde \S_2)<\infty)\le 2 P_y(T(\tilde \S_2)<T).
\ee
The argument is similar to the one showing $P_z(\tau_j)\le
2 P_z(\tau_j<T)$ in \reff{encage.11} 
of the proof of Lemma~\ref{lem-encageB}. We omit to reproduce it.
Thus, from \reff{trans.40} and \reff{trans.8},
\be{trans.41}
P_{y}(T(\S_2)<\infty)\le \frac{2^{d-1}c_1}{c_2} P_{y}(T(\tilde \S_2)<T).
\ee
We show that starting from $\tilde z\in \tilde \S_2$,
a walk has good chances of hitting $\tilde \S_1$ before $\Lambda$, 
as we show \reff{trans.40}, and here again we omit the argument showing
that for any $\tilde z_1\in \tilde \S_1$
\be{trans.17}
P_{\tilde z}(S_{T(\tilde \S_1)}=\tilde z_1)\le 2 P_{\tilde z}
\pare{T(\tilde \S_1)<T, S_{T(\tilde \S_1)}=\tilde z_1}.
\ee
\underline{ About the supremum in \reff{trans.3}}.
Now, by Harnack's inequality for the discrete Laplacian 
(see Theorem 1.7.2 of~\cite{LAWLER}),
there is $c_H>0$ independent of $n$ such that for any 
$z_2,z_2'\in \S_2$, and any $z_1\in \S_1$
\be{trans.16}
P_{z_2}\pare{S_{T(\S_1)}=z_1}\le c_H P_{z_2'}\pare{S_{T(\S_1)}=z_1}.
\ee
Now, using \reff{trans.17}, and the obvious fact
\[
P_{z_2'}\pare{S_{T(\S_1)}=z_1}=
P_{\T(z_2')}\pare{S_{T(\tilde \S_1)}=\T(z_1)}, 
\]
we obtain for any $z_1\in \S_1$
\be{trans.18}
\sup_{z\in\S_2} P_z(S_{T(\S_1)}=z_1)\le c_H
\inf_{\tilde z\in\tilde \S_2} P_{\tilde z}(S_{T(\tilde \S_1)}=\tilde z_1)
\le 2c_H \inf_{\tilde z\in\tilde \S_2} P_{\tilde z}
\pare{S_{T(\tilde \S_1)}=\tilde z_1,T(\tilde \S_1)<T}.
\ee
Starting with \reff{trans.3}, 
and combining \reff{trans.6}, \reff{trans.31},
and \reff{trans.18}, we obtain
\[
\begin{split}
P_y(S_T=x)& \le P_y(T(\S_2)<\infty)\sum_{z_1\in \S_1}
\pare{ \sup_{z\in \S_2}P_z(S_{T(\S_1)}=z_1)} P_{z_1}(S_T=x)\\
&\le \frac{2^{d-1}c_1}{c_2} P_{y}(T(\tilde \S_2)<T)
\sum_{\tilde z_1\in \tilde\S_1} \! 2c_H
\inf_{\tilde z\in\tilde \S_2} 
P_{\tilde z} \pare{S_{T(\tilde \S_1)}=\tilde z_1,T(\tilde \S_1)<T}\\
&\qquad \times C_0 P_{\tilde z_1}(S_T=\tilde x)\le P_y(S_T=\T(x)).
\end{split}
\]
\qed

%PROOF OF THE LEMMA.
\subsection{Proof of Lemma~\ref{improp-lem.1}.}
We only prove the first inequality in \reff{improp.21}, 
the second is similar. 
The proof uses arguments used in the proof of Lemma~\ref{lem-encageB}, and 
Lemma~\ref{trans-lem.1}. Namely, consider $x,x'\in \C$, and draw shells $\acc{C_k}$
and $\acc{D_k}$ as in \reff{encage.3} but around $x$ and $x'$ respectively. 
Note that here $C_k\cap D_{k'}$ may not be
empty. Also, choose $i$ and $j$ such that condition \reff{encage.4} holds.
Then, we decompose $\acc{S_T=x}$ by conditioning on $\S_1$ as in \reff{trans.3}.
On the term $P_{z_1}(S_T=x)$ we use the following  rough bound
\be{trans.30}
P_{z_1}(S_T=x)\le P_{z_1}(T(x)<\infty)\le \frac{c_d}{|z_1-x|^{d-2}}.
\ee
We now use the obvious observation that $2|z_1-x|\ge |z_1-x'|$.
Indeed, $|z_1-x|\ge \diam(\C)\ge |x-x'|$ implies that
$2|z_1-x|\ge |z_1-x|+|x-x'|\ge |z_1-x'|$ by the triangle inequality.
Thus there are a constant $c_3$ such that 
for the hitting time $\tau_j$ defined in \reff{encage.5}
\be{trans.31}
P_{z_1}(\tau_j<\infty)\ge \frac{c_2 \capa(D_j)}{|z_1-x'|^{d-2}}\ge
\frac{c_3}{|z_1-x|^{d-2}}.
\ee
From \reff{trans.3} and \reff{trans.31}, we have
\be{trans.32}
P_y(S_T=x)\le \frac{c_d}{c_3} \sum_{z_1\in \S_1} 
P_y\pare{T(\S_1)<T,S_{T(\S_1)}=z_1} P_{z_1}(\tau_j<\infty)
\ee
By argument \reff{encage.10}, and the choice of $l_0$ in \reff{encage.11}, we have
$2 P_{z_1}(\tau_j<T)\ge P_{z_1}(\tau_j<\infty)$. Finally, from $D_j$ to $x'$, there 
is a path avoiding $\Lambda'\bs\acc{x'}$ which cost a bounded amount depending only
on $l_0$.
\qed

%PROOF OF THE LEMMA.
\subsection{Proof of Corollary~\ref{improp-cor.1}.}
Note that by Lemma~\ref{improp-lem.1}, we have
\be{trans.33}
P_x(S_T=y)\le C_I P_x(S_T=y').
\ee
Now, $P_x(S_T=y')=P_{y'}(S_T=x)$, and we use again Lemma~\ref{improp-lem.1}
\be{trans.34}
P_{y'}(S_T=x)\le C_I P_{y'}(S_T=x')\Longrightarrow
P_x(S_T=y) \le C_I^2 P_{x'}(S_T=y').
\ee
\qed

\end{document}